\documentclass[11pt]{amsart}
\usepackage{palatino}
\usepackage{amsfonts}
\usepackage{amssymb}
\usepackage{amscd}

\theoremstyle{plain}

\newtheorem{theorem}{Theorem}[section]
\newtheorem{proposition}[theorem]{Proposition}
\newtheorem{corollary}[theorem]{Corollary}
\newtheorem{lemma}[theorem]{Lemma}
\theoremstyle{definition}

\newtheorem{remark}[theorem]{Remark}

\newtheorem{conjecture/question}[theorem]{Conjecture/Question}

\newtheorem{remark/definition}[theorem]{Remark/Definition}
\newtheorem{terminology/notation}[theorem]{Terminology/Notation}

\newcommand{\marginlabel}[1]%
  {\mbox{}\marginpar{\raggedleft\hspace{0pt}\bfseries\sf#1}}

\def\rmapdown#1{\Big\downarrow
   \rlap{$\vcenter{\hbox{$\scriptstyle#1$}}$ }}

\def\GG{{\textbf G}}

\def\PP{{\textbf P}}
\def\OO{\mathcal{O}}

\def\F{\mathcal{F}}

\def\E{\mathcal{E}}

\def\cM{\mathcal{M}}

\def\cC{\mathcal{C}}
\def\H{\mathcal{H}}
\def\Pic0{{\rm Pic}^0(X)}

\def\mm{\overline{\mathcal{M}}}
\def\cc{\overline{\mathcal{C}}}

\theoremstyle{remark}

\begin{document}

\title{Higher ramification and varieties of secant divisors on the generic curve}

\author[G. Farkas]{Gavril Farkas}

\address{Department of Mathematics, University of Texas,
Austin, TX 78712} \email{{\tt gfarkas@math.utexas.edu}}
\thanks{Research  partially supported by an Alfred P. Sloan Fellowship, the NSF Grants DMS-0450670 and DMS-0500747
and a 2006 Texas Summer Research Assignment}

\maketitle

\begin{center}
{\small{ABSTRACT}}
\end{center}

\begin{abstract}
For a smooth projective curve, the cycles of $e$-secant $k$-planes are among the most studied objects in classical enumerative geometry and there are well-known formulas due to Castelnuovo, Cayley and MacDonald concerning them. Despite various attempts, surprisingly little is known about the enumerative validity of such formulas. The aim of this paper is to  clarify this problem in the case of the generic curve $C$ of given genus. We determine precisely under which conditions the cycle of $e$-secant $k$-planes in non-empty and we compute its dimension. We also precisely determine the dimension of the variety of linear series on $C$ carrying $e$-secant $k$-planes.
\end{abstract}

\vskip 23pt

For a smooth projective curve $C$ of genus $g$, we denote by $C_e$
the $e$-th symmetric product of $C$ and by $G^r_d(C)$ the variety of
linear series of type $\mathfrak g^r_d$ on $C$, that is,
$$G^r_d(C):=\{(L, V): L\in \mbox{Pic}^d(C), V\in G(r+1, H^0(L))\}.$$
The  main result of Brill-Noether theory states that if $[C]\in
\cM_g$ is a general curve then $G^r_d(C)$ is a smooth variety of
dimension equal to $\rho(g, r, d):=g-(r+1)(g-d+r)$. For a linear
series $l=(L, V)\in G^r_d(C)$ and an effective divisor $D\in C_e$,
using the natural inclusion $H^0(C, L\otimes \OO_C(-D))\subset H^0(C, L)$,
we can define a new linear series $l(-D):=\bigl(L\otimes \OO_C(-D),
V\cap H^0(L\otimes \OO_C(-D))\bigr)$. We fix integers $0\leq f< e$
and introduce the determinantal cycle
$$V_e^{e-f}(l):=\{D\in C_e: \mbox{dim } l(-D)\geq r-e+f\}$$ of effective divisors of
degree $e$ which impose at most $e-f$ independent conditions on $l$.
If $l$ is very ample and we view $C\stackrel{l}\hookrightarrow
\PP^r$ as an embedded curve, then $V_e^{e-f}(l)$ parameterizes
$e$-secant $(e-f-1)$-planes to $C$. Each irreducible component of
$V_e^{e-f}(l)$ has dimension at least $e-f(r+1-e+f)$. The cycles
$V_e^{e-f}(l)$ have been extensively studied in classical
enumerative geometry. The virtual class $[V_e^{e-f}(l)]^{{virt}}\in
A^{f(r+1-e+f)}(C_e)$ has been computed by MacDonald and its
expression is tremendously complicated and thus of limited practical
use (see \cite{ACGH}, Chapter VIII). One case when we have a
manageable formula is for $e=2r-2$ and $f=r-1$, when
$[V_{2r-2}^{r-1}(l)]^{virt}$ computes the (virtual) number of
$(r-2)$-planes in $\PP^r$ which are $(2r-2)$-secant to $C$ (cf.
\cite{Ca}).

Surprisingly little is known about the validity of these classical
enumerative formulas (see \cite{H} and \cite{LB} for partial results
in the case of curves in $\PP^3$). The aim of this paper is to
clarify this problem for a general curve $[C]\in \cM_g$. For every
linear series $l\in G^r_d(C)$ we determine precisely under which
conditions the cycle $V_e^{e-f}(l)$ is non-empty and has the
expected dimension. Then having fixed $[C]\in \cM_g$, we determine
the dimension of the family of linear series $l\in G^r_d(C)$ with an
$e$-secant $(e-f-1)$-plane. For our first result, we use
degeneration techniques together with a few facts about the ample
cone of the moduli space $\mm_{0, g}$ to prove the following:
\begin{theorem}\label{secant}
Let $[C]\in \cM_g$ be a general curve and we fix non-negative integers
$0\leq f< e$,\ $r$
and $d$, such that $r-e+f\geq 0$. Then
we have that
$$\mathrm{dim }\{l\in G^r_d(C): V_e^{e-f}(l)\neq \emptyset\}\leq \rho(g, r,
d)-f(r+1-e+f)+e.$$ In particular, if $\rho(g,r,d)-f(r+1-e+f)+e<0$,
then $V_e^{e-f}(l)=\emptyset$, for every \ $l\in G^r_d(C)$.
\end{theorem}
\noindent More precisely, in Section 2 we prove the following
dimensionality estimate
$$\mbox{dim}\bigl\{(D, l)\in C_e\times G^r_d(C): D\in V_e^{e-f}(l)\bigr\}\leq
\rho(g, r, d)-f(r+1-e+f)+e,$$ which obviously implies Theorem
\ref{secant}.  This result generalizes the Brill-Noether theorem.
Indeed, when $l=K_C$, then $V_e^{e-f}(K_C)=C_e^f:=\{D\in C_e:
h^0(\OO_C(D))\geq f+1\}$. Since the fibres of the Abel-Jacobi map
$C_e^f\rightarrow W_e^f(C)$ are at least $f$-dimensional, clearly
$G_e^f(C)\neq \emptyset$ implies that $\mbox{dim } C_e^f\geq f$. Our
result reads $G_e^f(C)=\emptyset$ when $\rho(g, f, e)<0$, which is
 the non-existence part of the classical Brill-Noether theorem
(cf. \cite{EH1}). More generally, we have the following result in the case $\rho(g, r, d)=0$:
\begin{corollary}
Suppose $\rho(g, r, d)=0$ and $e<f(r+1-e+f)$. Then for a general
curve $[C]\in \cM_g$ we have that $V_e^{e-f}(l)=\emptyset$ for every
$l\in G^r_d(C)$, that is, no linear series of type $\mathfrak g^r_d$
on $C$ has any $e$-secant $(e-f-1)$-planes.
\end{corollary}

An immediate consequence of Theorem \ref{secant} is a proof of the
following conjecture of Coppens and Martens (cf. \cite{CM2} Theorem
3.3.1, for a proof in the case $f=1$):

\begin{corollary}\label{copm}
Let $[C]\in \cM_g$ be a general curve and we fix integers $0\leq
f<e$, $d$ and $r$ such that $r-e+f\geq 0$. Let $l$ be a general
linear series of type $\mathfrak g^r_d$ belonging to an irreducible
component of $G^r_d(C)$. Assuming that $V_e^{e-f}(l)$ is not empty,
then $e-f(r+1-e+f)\geq 0$. Moreover $V_e^{e-f}(l)$ is
equidimensional and $\mathrm{dim } \ V_e^{e-f}(l)=e-f(r+1-e+f)$.
\end{corollary}

We note that when $f=1$, Theorem \ref{secant} concerns the higher order very ampleness of
linear series on a general curve. We recall that a linear series $l\in G^r_d(C)$ is said to
 be $(e-1)$-very ample if $\mbox{dim } l(-p_1-\cdots -p_e)=r-e$, for any choice of (not necessarily
 distinct) $e$ points $p_1, \ldots, p_e\in C$. Thus $0$-very ampleness is equivalent to generation
 by global sections and $1$-very
 ampleness reduces to the classical notion of very ampleness.

\begin{corollary}\label{eample}
Let $[C]\in \cM_g$ be a general curve and $e, r, d$ be non-negative integers such that
 $\rho(g, r, d)+2e-2-r<0$. Then every linear series $l\in G^r_d(C)$ is
$(e-1)$-very ample.
\end{corollary}

Theorem \ref{secant} does not address the issue of existence of
linear series with $e$-secant $(e-f-1)$-planes. We prove the
following existence result for secant planes corresponding to linear
series $\mathfrak g^r_d$ on an \emph{arbitrary} smooth curve of
genus $g$.

\begin{theorem}\label{existence}
Let $[C]\in \cM_g$ be a general smooth curve and we fix integers
$0\leq f< e\leq g$, $d$ and $r$, \ such that $f(r+1-e+f)\geq e$,\
$d\geq 2e-f-1$,\ $\mbox{  }g-d+r\geq 0$,
$$\rho(g, r, d)-f(r+1-e+f)+e\geq 0 \mbox{ and } \rho(g, r-e+f,
d-e)\geq 0.$$ Assume moreover that we are in one of the following
situations: $$ (i)\ \mbox{ } 2f\leq e-1, \ (ii)\ \mbox{ }e=2r-2
\mbox{ and } f=r-1, \ \mbox{ } (iii)\  \mbox{ }e<2(r+1-e+f), \
\mbox{ or}
$$
$$ \ (iv)\mbox{ } \ \rho(g, r, d)\geq f(r+1-e+f)-(g-d+r).$$ Then
there exists a linear series $l\in G^r_d(C)$ such that
$V_e^{e-f}(l)\neq \emptyset$. Moreover, one has that the following
dimensionality statement:
$$\mathrm{dim}\{(D, l)\in C_e\times G^r_d(C): D\in
V_e^{e-f}(l)\}=\rho(g, r, d)-f(r+1-e+f)+e.$$
\end{theorem}

The inequalities $\rho(g, r-e+f, d-e)\geq 0$ and $\rho(g, r,
d)+e-f(r+1-e+f)\geq 0$ are obvious necessary conditions for the
existence of $l\in G^r_d(C)$ with $V_{e}^{e-f}(l)\neq \emptyset$ on
a general curve $[C]\in \cM_g$. To give an example, an elliptic
quartic curve $C\subset \PP^3$ has no $3$-secant lines even though
$\rho(g, r, d)+e-f(r+1-e+f)>0$ (note that $e=3$ and $f=1$ in this
case). Theorem \ref{existence} is stated in the range
$f(r+1-e+f)\geq e$, corresponding to the case when linear series
$l\in G^r_d(C)$ with $V_e^{e-f}(l)\neq \emptyset$ are expected to be
special in the Brill-Noether cycle $G^r_d(C)$. It is clear though
that the methods of this paper can be applied to the case $e\geq
f(r+1-e+f)$ as well. In that range however, when one expects
existence of $e$-secant $(e-f-1)$-planes for every $l\in G^r_d(C)$,
there are nearly optimal existence results obtained by using
positivity for Chern classes of certain vector bundles in the style
of \cite{FL}: For every curve $[C]\in \cM_g$ and $l\in G^r_d(C)$,
assuming that $d\geq 2e-1$\ and \ $e-f(r+1-e+f)\geq r-e+f$, one
knows that $V_e^{e-f}(l)\neq \emptyset$ (cf. \cite{CM}, Theorem
1.2). For $l\in G^r_d(C)$ such that $g-d+r\leq 1$ (e.g. when $l$ is
non-special), if we keep the assumption \ $e-f(r+1-e+f)\geq 0$, it
is known that $V_e^{e-f}(l)\neq \emptyset$ if and only if $\rho(g,
r-e+f, d-e)\geq 0$ (cf. \cite{ACGH}, pg. 356). This appears to be
the only case when MacDonald's formula displays some positivity
features that can be used to derive existence results on
$V_e^{e-f}(l)$. In the case $l=K_C$, one recovers of course the
existence theorem from classical Brill-Noether theory. We finally
mention that Theorem \ref{existence} holds independent of the
assumptions $(i)-(iii)$, whenever a certain transversality condition
(\ref{assumption2}) concerning a general curve $[Y, p]\in \cM_{e,
1}$ is satisfied (see Section 3 for details). Theorem
\ref{existence} is then proved by verifying this condition
(\ref{assumption2}) in each of the cases $(i)-(iii)$.

We now specialize to the case when $e=f(r+1-e+f)$ which is covered
by Theorem \ref{existence}. One can write $r=(u-1)(f+1)$ and $e=uf$
for some $u\geq 1$, and we obtain the following result concerning
the classical problem of existence of $uf$-secant secant
$(uf-f-1)$-planes to curves in $\PP^r$:

\begin{corollary}\label{existence2}
Let $C$ be a smooth curve of genus $g$. We fix integers $d, u, f\geq
2$ \ and assume that the inequalities \ $g\geq uf, \  d\geq
2uf-f-1$, \ $\rho(g, uf+u-f-1, d)\geq 0$ and $\rho(g, u-1, d-uf)\geq
0$  hold. Then there exists an embedding \ $C\subset
\PP^{(u-1)(f+1)}$ with $\mathrm{deg}(C)=d$, such that $C$ has a \
$uf$-secant $(uf-f-1)$-plane. If moreover, $[C]\in \cM_g$ is general
in moduli, then the embedded curve $C\stackrel{l}\hookrightarrow
\PP^{(u-1)(f+1)}$ corresponding to a general linear series $l\in
G^{(u-1)(f+1)}_d(C)$ \ has only a finite number of \ $uf$-secant
$(uf-f-1)$-planes.
\end{corollary}

If $[C]\in \cM_g$ is suitably general we can prove that the
Cayley-Castelnuovo numbers predicting the number of $(2r-2)$-secant
$(r-2)$-planes of a curve in $C\subset \PP^r$ have a precise
enumerative meaning:

\begin{theorem}\label{cayley}
Let $[C]\in \cM_g$ be a general curve. We fix integers $d, r\geq 3$
such that $d\geq 3r-2$, \ $\rho(g, r, d)\geq \emptyset$ and
$\rho(g, 1, d-2r+2)\geq 0$. Then if $C\stackrel{l}\hookrightarrow
\PP^r$ is an embedding corresponding to a general linear series
$l\in G^r_d(C)$, then $C$ has only finitely many $(2r-2)$-secant
$(r-2)$-planes. Their number (counted with multiplicities) is
$$C(d, g, r)=\sum_{i=0}^{r-1} \frac{(-1)^i}{r-i}{d-r-i+1\choose
r-1-i }{d-r-i \choose r-1-i}{g\choose i}.$$
\end{theorem}
A modern proof of the formula for $C(d, g, r)$ is due to MacDonald
and appears in \cite{ACGH} Chapter VIII. The original formula is due
to Castelnuovo (cf. \cite{Ca}). When $r=3$, we recover  Cayley's
formula for the number of $4$-secant lines of a smooth space curve
$C\subset \PP^3$ of degree $d$ (cf. \cite{C}):
$$C(d, g, 3)=\frac{1}{12} (d-2)(d-3)^2(d-4)-\frac{g}{2}
(d^2-7d+13-g).$$

To make a historical remark, there have been various attempts to
rigorously justify the so-called {\emph{functional method}} that
Cayley (1863), Castelnuovo (1889) and Severi (1900) used to derive
their enumerative formulas and to determine their range of
applicability (see \cite{LB}, \cite{V}). For instance, Cayley's
formula is shown to hold for an arbitrary smooth curve in $\PP^3$,
provided that $C(d, g, 3)$ is defined as the degree of a certain
$0$-cycle $\mbox{Sec}_4(C)$ in $\GG(1, 3)$ (cf. \cite{LB2}). The
drawback of this approach is that it becomes very difficult to
determine when this newly defined invariant is really enumerative.
For instance Le Barz only shows that this happens for very special
curves in $\PP^3$ (rational curves and generic complete
intersections) and one of the aims of this paper is to establish the
validity of such formulas for curves that are general with respect
to moduli.

The second topic  we study concerns ramification points of powers of
linear series on curves. This question appeared first in a
particular case in \cite{F}. We recall that for a pointed curve $[C,
p]\in \cM_{g, 1}$ and a linear series $l=(L, V)\in G^r_d(C)$, the
\emph{vanishing sequence} of $l$ at $p$
$$a^l(p): a_0^l(p)<\ldots < a_r^l(p)\leq d$$ is obtained by ordering
the set $\{\mbox{ord}_p(\sigma)\}_{\sigma\in V}$. The \emph{weight}
of $p$ with respect to $l$ is defined as $w^l(p):=\sum_{i=0}^r
(a_i^l(p)-i)$. One says that $p$ is a \emph{ramification point} of
$l$ if $w^l(p)\geq 1$ and we denote by $R(l)$ the finite set of
ramification points of $l$.  If $[C, p]\in \cM_{g, 1}$ and
$\overline{\alpha}: 0\leq \alpha_0\leq \ldots \leq \alpha_r\leq d-r$
is a Schubert index of type $(r, d)$,  the cycle
$$G^r_d(C, p, \overline{\alpha}):=\{l\in G^r_d(C): a_i^l(p)\geq
\alpha_i+i \mbox{ for } i=0\ldots r\}$$ can be realized as a
generalized determinantal variety inside $G^r_d(C)$ having virtual
dimension $\rho(g, r, d, \overline{\alpha}):=\rho(g, r,
d)-\sum_{j=0}^r \alpha_j$. For a general pointed curve $[C, p]\in
\cM_{g, 1}$, it is known that the virtual dimension equals the
actual dimension, that is,
$$\mbox{dim } G^r_d(C, p, \overline{\alpha})=\rho(g, r, d,
\overline{\alpha}) \mbox{ }\mbox{ (cf. \cite{EH2} Theorem 1.1).}$$
We address the following question: suppose $l=(L, V)\in G^r_d(C)$ is
a linear series with a prescribed ramification sequence
$\overline{\alpha}$ at a fixed point $p\in C$. Is then $p$ a
ramification point of any of the powers $L^{\otimes n}$ for $n\geq
2$?  If so, can we describe the sequence $a^{L^{\otimes n}}(p)$? One
certainly expects that under suitable genericity assumptions on $C$
and $L$, the points in $\bigcup_{n\geq 1} R(L^{\otimes n})$ should
be uniformly distributed on $C$. For example, it is known that for
every $C$ and $L\in \mbox{Pic}^d(C)$, the set $\bigcup_{n\geq 1}
R(L^{\otimes n})$ is dense in $C$ with respect to the classical
topology (cf. \cite{N}). Silverman and Voloch showed that for any
$L\in \mbox{Pic}^d(C)$ there exist finitely many points $p\in C$
such that the set $\{n\geq 1: p\in R(L^{\otimes n})\}$ is infinite
(cf. \cite{SV}).

We prove that on a generic pointed curve $[C, p]$, a linear series
$(L, V)$ and its multiples $L^{\otimes n}$ share no ramification
points, that is $R(l)$ and $R(L^{\otimes n})$ are as transverse as
they can be expected to be and moreover, the vanishing sequence
$a^{L^{\otimes n}}(p)$ is close to being  minimal:
\begin{theorem}\label{powers}
We fix a general pointed curve $[C, p]\in \cM_{g, 1}$, integers $r,
d\geq 1, n\geq 3$  and a Schubert index $\overline{\alpha}:0\leq
\alpha_0\leq \ldots \leq \alpha_r\leq d-r$. We also set
$m:=[(n+1)/2]$. Then for every linear series $l=(L, V)\in G^r_d(C,
p, \overline{\alpha})$ and every positive integer $$a< nd-\rho(g, r,
d, \overline{\alpha})-g-\bigl[\frac{g}{m}\bigr],$$ we have that
$h^0(C, L^{\otimes n}(-ap))=h^0(C, L^{\otimes n})-a=nd+1-g-a$. In
other words, $a_i^{L^{\otimes n}}(p)=i$ for $0\leq i\leq a-1$.
\end{theorem}
In the case $n=2$, when we compare $R(l)$ and $R(L^{\otimes 2})$ our
results are sharper:
\begin{theorem}\label{double}
We fix a general pointed curve $[C, p]\in \cM_{g, 1}$, integers $r,
d\geq 1$ and a Schubert index $\overline{\alpha}: \alpha_0\leq
\ldots \leq \alpha_r\leq d-r$. Then for every $(L, V)\in G^r_d(C, p,
\overline{\alpha})$ and every positive integer
$$a < \mathrm{max}\{2d+2-2g-\rho(g, r, d, \overline{\alpha})+\bigl[\frac{g-1}{2}\bigr],
\ \mbox{ } 2d+2-2g-2\rho(g, r, d,
\overline{\alpha})+2\bigl[\frac{g}{3}\bigr]\},$$ we have that
$h^0(C, L^{\otimes 2}(-ap))=h^0(C, L^{\otimes 2})-a=2d+1-g-a$.
\end{theorem}

Comparing the bounds on $a$ given in  Theorems \ref{powers} and
\ref{double} with the obvious necessary condition $a\leq nd-g+1$
which comes from the Riemann-Roch theorem, we see  that our results
are essentially optimal for relatively small values of $\rho(g, r,
d, \overline{\alpha})$ when the linear series $(L, V)\in G^r_d(C, p,
\overline{\alpha})$ have a strong geometric characterization. On the
other hand, if for instance $\rho(g, r, d, \overline{\alpha})=g$,
then $L\in \mbox{Pic}^d(C)$ and $p\in C$ are arbitrary and one
cannot expect to prove a uniform result about the vanishing of
$H^1(C, L^{\otimes n}(-a p))$.

Theorems \ref{powers} and \ref{double} concern line bundles $L$ with
 prescribed ramification at a given point $p\in C$. Such bundles
are of course very special in $\mbox{Pic}^d(C)$. If instead, we try
to describe $\bigcup_{n\geq 1} R(L^{\otimes n})$ for a general line
bundle $L\in \mbox{Pic}^d(C)$, the answer turns out to be
particularly simple. We give a short proof of the following result:
\begin{theorem}\label{genlb} Let $C$ be a smooth curve of genus $g$ and $L\in
\rm{Pic}$$^d(C)$ a very general line bundle. \newline \noindent (1)
All the ramification points of the powers $L^{\otimes n}$ are
ordinary, that is, $w^{L^{\otimes n}}(p)\leq 1$ for all $p\in C$ and
$n\geq 1$. \newline \noindent (2) $R(L^{\otimes a})\cap R(L^{\otimes
b})=\emptyset$ for $a\neq b$, that is, a point $p\in C$ can be a
ramification point for at most a single power of $L$.
\end{theorem}
After this paper has been written I have learnt that Theorem
\ref{genlb} has also been proved independently by M. Coppens in
\cite{Co}. I would like to thank the referee for a very careful
reading of this paper and for pointing out that the initial proof of
Theorem \ref{secant} was not complete.

\section{Ramification points of multiples of linear series}
In this section we use the technique of limit linear series to prove
Theorems \ref{powers} and \ref{double}. We start by fixing a
Schubert index $\overline{\alpha}: 0\leq \alpha_0\leq \ldots \leq
\alpha_r\leq d-r$ and two integers $a\geq 0, n\geq 2$. We also set
$m:=[(n+1)/2]$.

We assume that for every $[C, p]\in \cM_{g, 1}$ there exists a
linear series $l=(L, V)\in G^r_d(C, p, \overline{\alpha})$ such that
$H^0(K_C\otimes L^{\otimes (-n)}\otimes \OO_C(ap))\neq 0$. By a
degeneration argument we are going to show that this implies the
inequalities
\begin{equation}\label{ineq1}
a\geq nd-g-\rho(g, r, d,
\overline{\alpha})-\bigl[\frac{g}{m}\bigr],\ \ \mbox{ when }n\geq 3,
\end{equation}
\begin{equation}\label{ineq2}
a\geq 2d+2-2g- \rho(g, r, d,
\overline{\alpha})+\bigl[\frac{g-1}{2}\bigr], \ \
\end{equation}
and \begin{equation} \label{numerology} a\geq 2d+2-2g-2\rho(g, r, d,
\overline{\alpha})+2\bigl[\frac{g}{3}\bigr], \ \ \mbox{when } n=2.
\end{equation} This will prove both Theorems \ref{powers} and
\ref{double}.

We degenerate $[C, p]$ to a stable curve $[X_0:=E_0\cup_{p_1} E_1
\cup_{p_2} \ldots\cup_{p_{g-1}} E_{g-1}, p_0]$, where $E_i$ is a
general elliptic curve, $p_i, p_{i+1}\in E_i$ are points such that
$p_{i+1}-p_i\in \mbox{Pic}^0(E_i)$ is not a torsion class  and
moreover $E_i\cap E_{i+1}=\{p_{i+1}\}$ for $0\leq i\leq g-2$. Thus
$X_0$ is a string of $g$ elliptic curves and the marked point $p_0$
specializes to a general point lying on the first component $E_0$.
We also consider a $1$-dimensional family
$\pi:\mathcal{X}\rightarrow B$ together with a section
$\sigma:B\rightarrow \mathcal{X}$, such that $B=\mbox{Spec}(R)$ with
$R$ being a discrete valuation ring having uniformizing parameter
$t$. We assume that $\mathcal{X}$ is a smooth surface and that there
exists an isomorphism between $X_0$ and $\pi^{-1}(0)$. Under this
isomorphism we also assume that $\sigma(0)=p_0\in X_0$.  Here $0\in
B$ is the point corresponding to the maximal ideal of $R$ and we
denote by $\eta$ and $\overline{\eta}$ the generic and geometric
generic point of $B$ respectively. By assumption, there exists a
linear series $l_{\overline{\eta}}=(L_{\overline{\eta}},
V_{\overline{\eta}})\in G^r_d(X_{\overline{\eta}},
\sigma(\overline{\eta}), \overline{\alpha}),$ such that
$H^0(X_{\overline{\eta}}, \omega_{X_{\overline{\eta}}}\otimes
L_{X_{\overline{\eta}}}^{\otimes (-n)}\otimes
\OO_{X_{\overline{\eta}}}(a\sigma(\overline{\eta})))\neq 0$. By
possibly blowing up $\mathcal{X}$ at the nodes of $X_0$ and thus
replacing the central fibre by a curve $X$ obtained from $X_0$ by
inserting chains of smooth rational curves at the points $p_1,
\ldots, p_{g-1}$, we may assume that $l_{\overline{\eta}}$ comes
from a linear series $l_{\eta}=(L_{\eta}, V_{\eta})\in
G^r_d(X_{\eta}, \sigma(\eta), \overline{\alpha})$ on the generic
fibre $X_{\eta}$.

We denote by $l_{E_i}=(L_{E_i}, V_{E_i})\in G^r_d(E_i)$ the
$E_i$-aspect of the limit linear series on $X$ induced by
$l_{\eta}$: Precisely,  if $\mathcal{L}$ is a line bundle on
$\mathcal{X}$ extending $L_{\eta}$, then $L_{E_i}\in
\mbox{Pic}^d(E_i)$ is the restriction to $E_i$ of the unique twist
$\mathcal{L}_{E_i}$ of $\mathcal{L}$ along components of
$\pi^{-1}(0)$ such  that $\mbox{deg}_Z(\mathcal{L}_{i |Z})=0$ for
any irreducible component $Z\neq E_i$ of $\pi^{-1}(0)$ (see also
\cite{EH1}, p. 348). Since we gave ourselves the freedom of
blowing-up $\mathcal{X}$ at the nodes of $\pi^{-1}(0)$, we can also
assume that $\{l_{E_i}\}_{i=0}^{g-1}$ constitutes a limit $\mathfrak
g^r_d$ on $X_0$ which is obtained from a refined limit $\mathfrak
g^r_d$ on $X$ by retaining only the aspects of the elliptic
components of $X$. The compatibility relations between the vanishing
orders of the $l_{E_i}$'s imply the following inequality between
Brill-Noether numbers:
\begin{equation}\label{bn}
\rho(g, r, d, \overline{\alpha})\geq \rho(l_{E_0}, p_0,
p_1)+\rho(l_{E_1}, p_1, p_2)+\cdots +\rho(l_{E_{g-2}}, p_{g-2},
p_{g-1})+\rho(l_{E_{g-1}}, p_{g-1}),
\end{equation}
where $\rho(l_{E_i}, p_i, p_{i+1}):=\rho(1, r,
d)-w^{l_{E_i}}(p_i)-w^{l_{E_{i}}}(p_{i+1})$. By assumption, there
exists a non-zero section $\rho_{\eta}\in H^0\bigl(X_{\eta},
\omega_{X_{\eta}}\otimes \mathcal{L}_{\eta}^{\otimes (-n)}\otimes
\OO_{X_{\eta}}(a \sigma(\eta))\bigr)$. This implies that  if we
denote by $\tilde{\mathcal{L}}_i$ the unique line bundle on the
surface $\mathcal{X}$ such that  (1) $\tilde{\mathcal{L}}_{i |
X_{\eta}}=L_{\eta}$, and  (2) $\mbox{deg}_Z(\omega_X\otimes
\tilde{\mathcal{L}_i}^{\otimes (-n)}\otimes \OO_X(ap_0))=0$, for
every component $Z$ of $X$ such that $Z\neq E_i$, then $H^0(E_i,
\omega_{X}\otimes \tilde{\mathcal{L}}_{i}^{\otimes (-n)}\otimes
\OO_X(a p_0)\otimes \OO_{E_i})\neq 0$. We set
$$\cM_i:=\omega_{\pi}\otimes \tilde{\mathcal{L}}_i^{\otimes
(-n)}\otimes \OO_{\mathcal{X}}(a\ \sigma(B))\in
\mbox{Pic}(\mathcal{X}).$$ Then $\cM_{i
|E_i}=\OO_{E_i}\bigl((a+2i)\cdot p_i+(2g-2-2i)\cdot p_{i+1}\otimes
L_{E_i}^{\otimes (-n)}\bigr)$ for all $0\leq i\leq g-1$. For each
such $i$ we denote by $n_i$ the smallest integer such that
$\tilde{\rho}_i:=t^{n_i}\rho_{\eta} \in \pi_*(\cM_i)$ and we set
$$\rho_i:=\tilde{\rho}_{i |E_i} \in H^0(E_i, \cM_{i |E_i}).$$ Thus
$0\neq \rho_i \in H^0(E_i, \OO_{E_i}((a+2i)\cdot p_i+(2g-2-2i)\cdot
p_{i+1}\otimes L_{E_i}^{\otimes (-n)}))$ and in a way similar to
\cite{EH1} Proposition 2.2, we can prove that
\begin{equation}\label{limitlin}
\mbox{ord}_{p_i}(\rho_i)+\mbox{ord}{p_i}(\rho_{i-1})\geq
2g-2-nd+a=\mbox{deg}(\cM_{i |E_i}).
\end{equation}
One also has the inequalities
$\mbox{ord}_{p_i}(\rho_i)+\mbox{ord}_{p_{i+1}}(\rho_i)\leq
2g-2-nd+a$ (and similar inequalities when passing through the
rational components of $X$), from which it follows that one can
write down a non-decreasing sequence of vanishing orders
\begin{equation}\label{sequence}
0\leq \mbox{ord}_{p_0}(\rho_0)\leq \mbox{ord}_{p_1}(\rho_1)\leq
\ldots \leq \mbox{ord}_{p_i}(\rho_i)\leq \ldots \leq
\mbox{ord}_{p_{g-1}}(\rho_{g-1}).
\end{equation}
 Since $\rho_{g-1}$
is a non-zero section of a line bundle of degree $2g-2-nd+a$ on $E_{g-1}$, we
must have that $\mbox{ord}_{p_{g-1}}(\rho_{g-1})\leq 2g-2-nd+a$.
This inequality will eventually lead to the bound on the constant
$a$.

Let us suppose now that we have fixed one of the elliptic components of
$X$, say $E_i$, such that $\rho(l_{E_i}, p_i, p_{i+1})=0$. By
counting dimensions, we see that for every $0\leq j\leq r$ there
exists a section $u_j\in V_{E_i}$ such that $\mbox{div}(u_j)\geq
a_j^{l_{E_i}}(p_i)\cdot p_i+a_{r-j}^{l_{E_i}}(p_{i+1})\cdot
p_{i+1}$. In particular, we have that
$a_j^{l_{E_i}}(p_i)+a_{r-j}^{l_{E_i}}(p_{i+1})\leq d$. Since
$p_{i+1}-p_i\in \mbox{Pic}^0(E_i)$ is not a torsion class, it
follows that the equality
$a_j^{l_{E_i}}(p_i)+a_{r-j}^{l_{E_i}}(p_{i+1})=d$ can hold for at
most one value $0\leq j\leq r$. Because $\rho(l_{E_i}, p_i,
p_{i+1})=0$, this implies that
$$a_{j}^{l_{E_i}}(p_i)+a_{r-j}^{l_{E_i}}(p_{i+1})\geq d-1 \mbox{  for all
 } \ 0\leq j\leq r,$$ and there exists precisely one such index $j$ such
that $a_j^{l_{E_i}}(p_i)+a_{r-j}^{l_{E_i}}(p_{i+1})=d$. In this case
we get that $\mbox{div}(u_j)=a_j^{l_{E_i}}(p_i)\cdot
p_i+a_{r-j}^{l_{E_i}}(p_{i+1})\cdot p_{i+1}$, and for degree reasons
we must have that $L_{E_i}=\OO_{E_i}(a_j^{l_{E_i}}(p_i)\cdot
p_i+a_{r-j}^{l_{E_i}}(p_{i+1})\cdot p_{i+1})\in \mbox{Pic}^d(E_i).$

To summarize, if $\rho(l_{E_i}, p_i, p_{i+1})=0$, then the vanishing
sequence $a^{l_{E_{i+1}}}(p_{i+1})$ of the $E_{i+1}$-aspect of the
limit $\mathfrak g^r_d$ on $X$, is obtained from the vanishing
sequence $a^{l_{E_i}}(p_i)$ by raising all entries by $1$, except
one single entry which remains unchanged. Thus,
$a_j^{l_{E_i}}(p_i)=a_j^{l_{E_{i+1}}}(p_{i+1})$ for one index $0\leq
j\leq r$ and $a_k^{l_{E_{i+1}}}(p_{i+1})=a_k^{l_{E_i}}(p_i)+1$ for
$k\neq j$.

We now study what happens to the non-decreasing sequence
(\ref{sequence}) as we pass through a component $E_i$ with
$\rho(l_{E_i}, p_i, p_{i+1})=0$. Assume that
$\mbox{ord}_{p_i}(\rho_i)=\mbox{ord}_{p_{i+1}}(\rho_{i+1}):=b$. This
implies that $\mbox{ord}_{p_{i+1}}(\rho_i)=2g-2-nd+a-b$ and
$$L_{E_i}^{\otimes n}=\OO_{E_i}((a+2i-b)\cdot p_i+(nd-a+b-2i)\cdot
p_{i+1})\in \mbox{Pic}^{nd}(E_i).$$ Because $\rho(l_{E_i}, p_i,
p_{i+1})=0$, as we have seen, $L_{E_i}$ can be represented by an
effective divisor which is supported only at $p_i$ and $p_{i+1}$.
Precisely, we can write that
$L_{E_i}=\OO_{E_i}\bigl(a_j^{l_{E_i}}(p_i)\cdot
p_i+a_{r-j}^{l_{E_i}}(p_{i+1})\cdot p_{i+1}\bigr)$ for a unique
$0\leq j\leq r$. Since $L_{E_i}$ cannot admit two different
representations by effective divisors supported only at $p_i$ and
$p_{i+1}$, we must have that
\begin{equation}\label{divisible}
L_{E_i}=\OO_{E_i}\Bigl(\frac{a+2i-b}{n}\cdot
p_i+\frac{nd-a+b-2i}{n}\cdot p_{i+1}\Bigr). \end{equation}
 In
particular, we have that $(a+2i-b)/n \in \mathbb Z$ and
$a_j^{l_{E_i}}(p_i)=(a+2i-b)/n $.

We consider a connected subcurve $Y\subset X$ containing $m+1$
elliptic components $E_i$ and we measure the increase in
(\ref{sequence}) as we pass through the components of $Y$.

\begin{lemma}\label{inec}
We fix $m:=[(n+1)/2]$ and integers $i$ and $b$ such that $bm\leq
i\leq g-1$. We denote by $R(i):=\#\{0\leq l\leq i-1: \rho(l_{E_l},
p_{l}, p_{l+1})\geq 1\}$. Then the following inequality holds:
$$\mathrm{ord}_{p_i}(\rho_i)+R(i)\geq b(m-1).
$$
\end{lemma}
\begin{proof}
We proceed by induction on $b$. For $b=0$ there is nothing to prove.
We set $b\geq 1$, $i:=(b-1)m$ and we assume that
$\mbox{ord}_{p_i}(\rho_i)+R(i)\geq (b-1)(m-1)$. We are going to
prove that the following inequality holds:
\begin{equation}\label{ind}
\mathrm{ord}_{p_{i+m}}(\rho_{i+m})-\mathrm{ord}_{p_i}(\rho_i)+R(i+m)-R(i)\geq
m-1.
\end{equation}

Assume this is not the case. Then there exist integers $0\leq l<
j\leq m-1$ such that the following relations hold: \ $(i)\
\rho(l_{E_{i+l}}, p_{i+l}, p_{i+l+1})=\rho(l_{E_{i+j}}, p_{i+j},
p_{i+j+1})=0$ \ and
$$(ii)\
\mbox{ord}_{p_{i+l}}(\rho_{i+l})=\mbox{ord}_{p_{i+l+1}}(\rho_{i+l+1}):=b,
\mbox { } \
\mbox{ord}_{p_{i+j}}(\rho_{i+j})=\mbox{ord}_{p_{i+j+1}}(\rho_{i+j+1}):=c.$$
Using (\ref{divisible}) this implies that $$
L_{E_{i+l}}=\OO_{E_{i+l}}\bigl(\frac{a+2i+2l-b}{n}\cdot
p_{i+l}+\frac{nd-a+b-2i-2l}{n}\cdot p_{i+l+1}\bigr), \ \mbox{ and
}$$
$$L_{E_{i+j}}=\OO_{E_{i+j}}\bigl(\frac{a+2i+2j-c}{n}\cdot
p_{i+j}+\frac{nd-a+c-2i-2j}{n}\cdot p_{i+j+1}\bigr).$$ In
particular, $(2j-2l-c+b)/n \in \mathbb Z$, hence we can write
$c=b-kn+2(j-l)$ for some $k\in \mathbb Z$. If $k\geq 1$, since
$c\geq b$, we obtain that $m-1\geq j-l \geq n/2$, which is a
contradiction. Therefore we must have that $k\leq 0$, and this holds
for every pair $(j, l)$ satisfying (i) and (ii). We choose now the
pair $0\leq l<j\leq m-1$ satisfying (i) and (ii) and for which
moreover, the difference $j-l$ is maximal.

For each integer $0\leq e\leq l-1$ we have that either
$\rho(l_{E_{i+e}}, p_{i+e}, p_{i+e+1})\geq 1$ or
$\mbox{ord}_{p_{i+e+1}}(\rho_{i+e+1})>\mbox{ord}_{p_{i+e}}(\rho_{i+e})$.
This fact leads to the inequality
\begin{equation}\label{num1}
\mbox{ord}_{p_{i+l}}(\rho_{i+l})-\mbox{ord}_{p_i}(\rho_i)+R(i+l)-R(i)\geq
l.
\end{equation} Similarly, by studying the subcurve of $Y$ containing
$E_{i+j+1}, \ldots, E_{i+m-1}$, we find that
\begin{equation}\label{num2}
\mbox{ord}_{p_{i+m}}(\rho_{i+m})-\mbox{ord}_{p_{i+j+1}}(\rho_{i+j+1})+R(i+m)-R(i+j+1)\geq
m-j-1.\end{equation} Finally, we look at the subcurve of $X$
containing $E_{i+l}, \ldots, E_{i+j}$ and we can write
\begin{equation}\label{num3}
\mbox{ord}_{p_{i+j}}(\rho_{i+j})-\mbox{ord}_{p_{i+l}}(\rho_{i+l})+R(i+j+1)-R(i+l)\geq
c-b\geq 2(j-l)\geq j-l+1.
\end{equation} By adding (\ref{num1}), (\ref{num2}) and (\ref{num3})
together we obtain (\ref{ind}) which proves the Lemma.
\end{proof}

When $n=2$ we have a slightly better estimate than in the general
case:

\begin{lemma}\label{inec2}
($n=2$)   (1) Let $i$ be an integer such that $2b\leq i\leq g-1$.
Then $\mathrm{ord}_{p_i}(\rho_i)+R(i)\geq b$. \newline \noindent (2)
We fix $0\leq i\leq g-4$ and let $Y$ be a connected subcurve of $X$
containing precisely three elliptic curves $E_i, E_{i+1}$ and
$E_{i+2}$. If $R(i+3)=R(i)$, that is,
$$\rho(l_{E_i}, p_i, p_{i+1})=\rho(l_{E_{i+1}}, p_{i+1},
p_{i+2})=\rho(l_{E_{i+2}}, p_{i+2}, p_{i+3})=0,$$ then we have the
inequality $\mathrm{ord}_{p_{i+3}}(\rho_{i+3})\geq
\mathrm{ord}_{p_i}(\rho_i)+2.$
\end{lemma}
\begin{proof}
We only prove (2), the remaining statement being analogous to Lemma
\ref{inec}. We may assume that
$\mbox{ord}_{p_i}(\rho_i)=\mbox{ord}_{p_{i+1}}(\rho_{i+1}):=b$.
Hence $(a+2i-b)/2 \in \mathbb Z$ and there exists an index $0\leq
j\leq r$ such that
$$a_j^{l_{E_i}}(p_i)=a_j^{l_{E_{i+1}}}(p_{i+1})=\frac{a+2i-b}{2}, \mbox{ while
 }
\ a_k^{l_{E_{i+1}}}(p_{i+1})=a_k^{l_{E_i}}(p_i)+1 \ \mbox{ for }
k\neq j.$$ If
$\mbox{ord}_{p_{i+2}}(\rho_{i+2})=\mbox{ord}_{p_{i+1}}(\rho_{i+1})=b$,
then (\ref{divisible}) implies that $(a+2i+2-b)/2$ is an entry in
the vanishing sequence $a^{l_{E_{i+1}}}(p_{i+1})$. But this is
impossible, because $(a+2i-b)/2$ was an entry in the sequence
$a^{l_{E_i}}(p_i)$, hence we must have that
$\mbox{ord}_{p_{i+2}}(\rho_{i+2})\geq b+1$. Next, if
$\mbox{ord}_{p_{i+3}}(\rho_{i+3})=b+1$, this implies that
$\mbox{ord}_{p_{i+3}}(\rho_{i+3})=\mbox{ord}_{p_{i+2}}(\rho_{i+2})=b+1$,
hence again $\bigl(a+2(i+2)-(b+1)\bigr)/2 \in \mathbb Z$, which is
not possible for parity reasons. Thus we must have that
$\mbox{ord}_{p_{i+3}}(\rho_{i+3})\geq b+2$.
\end{proof}

\noindent {\emph{Proof of Theorem \ref{powers}.}} We complete the
proof of our result in the case $n\geq 3$. We write $g=bm+c$ with
$0\leq c\leq m-1$ and we set $i:=bm$. From Lemma \ref{inec} we
obtain that $\mbox{ord}_{p_i}(\rho_i)+R(i)\geq b(m-1)$.  Using the
reasoning of Lemma \ref{inec} for the connected subcurve of $X$
which contains $E_i, E_{i+1}, \ldots, E_{i+c-1}=E_{g-1}$, we get
that
\begin{equation}\label{num4}
\mbox{ord}_{p_{g-1}}(\rho_{g-1})-\mbox{ord}_{p_i}(\rho_i)+R(g-1)-R(i)\geq
c-2.
\end{equation}
Using (\ref{num4}), together with the inequality $R(g-1)\leq \rho(g,
r, d, \overline{\alpha})$, we can write that
$$\mbox{deg}(K_C\otimes L^{\otimes (-n)}\otimes \OO_C(ap))=
2g-2-nd+a\geq \mbox{ord}_{p_{g-1}}(\rho_{g-1})\geq
g-\bigl[\frac{g}{m}\bigr]-\rho(g, r, d, \overline{\alpha})-2,$$
which finishes the proof of Theorem \ref{powers}. \hfill $\Box$

\noindent{\emph{Proof of Theorem \ref{double}}}. From Lemma
\ref{inec2}  part (1), we obtain that
$$\mbox{ord}_{p_{g-1}}(\rho_{g-1})+R(g-1)\geq [(g-1)/2].$$ Since
$R(g-1)\leq \rho(g, r, d, \overline{\alpha})$, this leads to the
inequality $a\geq 2d+2-2g+[(g-1)/2]-\rho(g, r, d,
\overline{\alpha})$. To prove (\ref{numerology}) we divide $X$ into
$e:=[g/3]+1$ connected subcurves $Y_1, \ldots, Y_{e}$ such that
$Y_1, \ldots, Y_{e-1}$ each contain three elliptic components,
$\#(Y_i\cap Y_{i+1})=1$ for all $1\leq i\leq e-2$ and
$Y_e:=\overline{(\cup_{i=1}^{e-1} Y_i)^c}$. The curves $Y_i$ fall
into two categories: those for which there exists an elliptic
component $E_l\subset Y_i$ such that $\rho(l_{E_l}, p_l,
p_{l+1})\geq 1$ (and there are at most $\rho(g, r, d,
\overline{\alpha})$ such $Y_i$'s), and those for which
$\rho(l_{E_l}, p_l, p_{l+1})=0$ for each elliptic component
$E_l\subset Y_i$. Lemma \ref{inec2} part (2) gives that
$\mbox{ord}_{p_{g-1}}(\rho_{g-1}) \geq 2([g/3]-\rho(g, r, d,
\overline{\alpha}))$. This proves (\ref{ineq2}) and finishes the
proof of Theorem \ref{double}. \hfill $\Box$
\begin{remark} It is natural to ask how close to being optimal are the
bounds we obtained above. For $\rho(g, r, d, \overline{\alpha})$
relatively small, when any $L\in G^r_d(C, p, \overline{\alpha})$ has
a strong geometric characterization, the inequalities (\ref{ineq1}),
(\ref{ineq2}) and (\ref{numerology}) are in fact optimal. To see an
example, we set $g=3, r=3, d=6$ and $\rho(g, r, d,
\overline{\alpha})=0$. Thus we look at $\mathfrak g^3_6$'s on a
general $[C, p]\in \cM_{3, 1}$ having ramification at $p$ equal to
$(0\leq \alpha_0\leq \alpha_1\leq \alpha_2\leq \alpha_3\leq 3)$,
where $\sum_{i=0}^3 \alpha_i=3$. Theorem \ref{double} gives us that
$H^0(K_C\otimes L^{\otimes (-2)}\otimes \OO_C(a\cdot p))=0$ for
every integer $a\leq 9$. We show that this is optimal by noting that
when $a=10$ and $\overline{\alpha}=(0, 0, 1, 2)$, we have that
$$H^0(K_C\otimes L^{\otimes (-2)}\otimes \OO_C(10p))\neq 0, \mbox{
for every }L\in G^3_6(C, p, \overline{\alpha}).$$ Indeed, any such
linear series is of the form $L=K_C\otimes A^{\vee}\otimes \OO_C(5
p)\in W^3_6(C)$, where $A\in W^1_3(C)$ is such that $h^0(A(-2p))\geq
1$. A non-hyperelliptic curve of genus $3$ has two such $\mathfrak
g^1_3$'s. Precisely, if $z, t\in C$ are the two points the tangent
line at $p$ to $C\stackrel{|K_C|} \hookrightarrow \PP^2$ meets $C$
again, then $A=\OO_C(2p+z)$ or $A=\OO_C(2p+t)$. Say, we choose
$A=\OO_C(2p+z)$. By direct calculation we obtain that $ L^{\otimes
2}\otimes \OO_C(-10p)=K_C^{\otimes 2}\otimes A^{\otimes
(-2)}=\OO_C(2t)$, hence $h^0(K_C\otimes L^{\otimes (-2)}\otimes
\OO_C(10 p))=1$.
\end{remark}

\section{Varieties of secant planes to the general curve}
We fix a smooth curve $[C]\in \cM_g$ and two integers $0\leq f<e$.
In this section we study the varieties $V_e^{e-f}(l)$ of $e$-secant
$(e-f-1)$-planes corresponding to a linear series $l\in G^r_d(C)$.
We first define the correspondence $$\Sigma_C:=\{(D, l)\in C_e\times
G^r_d(C): \mbox{dim } l(-D)\geq r-e+f\},$$ and denote by
$\pi_1:\Sigma_C\rightarrow C_e$ and $\pi_2:\Sigma_C\rightarrow
G^r_d(C)$ the two projections. We assume that $\Sigma_C\neq
\emptyset$ for the general curve $[C]\in \cM_g$. Under this
assumption, we show that
\begin{equation}\label{sec1}
\mbox{dim}(\Sigma_C)\leq \rho(g, r,d)-f(r+1-e+f)+e.
\end{equation}
(We recall that the dimension of a scheme is the maximum of the
dimensions of its irreducible components). Since $\Sigma_C$ is a
determinantal subvariety of $C_e\times G^r_d(C)$, it follows that
for a general $[C]\in \cM_g$, if non-empty, the scheme $\Sigma_C$ is
equidimensional and $\mbox{dim}(\Sigma_C)=\rho(g, r,
d)-f(r+1-e+f)+e$. Note that this result does not establish the
non-emptiness of $\Sigma_C$ which is an issue that we will deal with
in Section 3. In any event, (\ref{sec1}) implies the dimensional
estimate
$$\mbox{dim}\{l\in G^r_d(C): V^{e-f}_e(l)\neq \emptyset\}\leq
\rho(g, r, d)-f(r+1-e+f)+e. $$ This will prove Theorem \ref{secant}
as well as Corollaries \ref{copm} and \ref{eample}.

We start by setting some notation. We denote by $j:\mm_{0,
g}\rightarrow \mm_g$ the ``flag'' map obtaining by attaching to each
stable curve $[R, x_1, \ldots, x_g]\in \mm_{0, g}$ fixed elliptic
tails $E_1, \ldots, E_g$ at the points $x_1, \ldots, x_g$
respectively. Thus $j([R, x_1, \ldots,
x_g]):=[\tilde{R}]=[R\cup_{x_1} E_1\cup \ldots \cup_{x_g} E_g]$ and
for such a curve, we denote by $p_R:\tilde{R}\rightarrow R$ the
projection onto $R$, that is, $p_R(E_i)=\{x_i\}$ for $1\leq i\leq
g$. We denote by $\cc_{g, n}=\mm_{g,n+1}$ the universal curve and by
$\pi:\cc_{g, n}\rightarrow \mm_{g, n}$ the morphism forgetting the
$(n+1)$-st marked point. We  write $\pi_e:\cc_{g,n}^e \rightarrow
\mm_{g, n}$ for the $e$-fold fibre product of $\cc_{g, n}$ over
$\mm_{g, n}$ and we introduce a map $\chi:\mm_{0, g}\times_{\mm_g}
\cc_g^e\rightarrow \cc_{0, g}^e$ which collapses the elliptic tails.
Thus $\chi$ is defined by
$$\chi\bigl([R, x_1,\ldots, x_g], (y_1, \ldots, y_e)\bigr):=\bigl([R, x_1,
\ldots, x_g], p_R(y_1), \ldots, p_R(y_e)\bigr),$$ for points $y_1, \ldots,
y_e \in \tilde{R}$. Let $W\subset \cc_g^e$ be the closure of the
locus
$$\{[C, y_1, \ldots, y_e]\in \cC_{g}^e: \exists l\in G^r_d(C)\mbox{
with } \mbox{ dim } l(-y_1-\cdots -y_e)\geq r-e+f\}.$$ By assumption
$\pi_e(W)=\mm_g$ and we define the locus $U:=\chi\bigl(W\cap
(\mm_{0, g}\times _{\mm_g}\cc_g^e)\bigr)$. Then $\pi_e(U)=\mm_{0,
g}$ and we denote by $e-m$ the minimal fibre dimension of the map
$\pi_{e | U}: U\rightarrow \mm_{0, g}$. Thus $0\leq m\leq e$ and
$\mbox{dim}(U\cap \pi_e^{-1}[R, x_1, \ldots, x_g])\geq e-m$, for
every $[R, x_1, \ldots, x_g]$, with equality for a general point
$[R, x_1, \ldots, x_g]\in \mm_{0, g}$.

We recall that for every choice of $4$ marked points $\{i, j, k,
l\}\subset \{1, \ldots, g\}$, one has a fibration $\pi_{i j k
l}:\mm_{0, g}\rightarrow \mm_{0, 4}$ obtained by forgetting the
marked points with labels in the set $\{i, j, k, l\}^c$ and
stabilizing the resulting rational curve. If we single out the first
$3$ marked points $x_1, x_2, x_3$ as being $0, 1$ and $\infty$, in
this way we obtain a birational map $\pi_{123}=(\pi_{1234},\ldots,
\pi_{123i}, \ldots, \pi_{123g}):\mm_{0, g}\rightarrow \mm_{0,
4}^{g-3}=(\PP^1)^{g-3}$ defined by
$$\pi_{123}([R, x_1, \ldots, x_g]):=\bigl([R, x_1, x_2, x_3, x_4],
[R, x_1, x_2, x_3, x_5], \ldots, [R, x_1, x_2, x_3, x_g]\bigr).$$
The map $\pi_{123}$ expresses $\mm_{0, g}$ as a blow-up of
$(\PP^1)^{g-3}$ such that all exceptional divisors of $\pi_{123}$
are boundary divisors of $\mm_{0, g}$ (cf. \cite{K}). In a similar
manner, one has a birational map $f:\cc_{0, g}^e\rightarrow \mm_{0,
4}^{g-3+e}=(\PP^1)^{g-3+e}$ defined by $f\bigl([R, x_1, \ldots,
x_g], y_1, \ldots, y_e\bigr):=$ $$:=\bigl([R, x_1, x_2, x_3, x_4],
\ldots, [R, x_1, x_2, x_3, x_g], [R, x_1, x_2, x_3, y_1], \ldots,
[R, x_1, x_2, x_3, y_e]\bigr).$$ For simplicity, sometimes we write
$f([R, x_1, \ldots, x_g], y_1, \ldots, y_e)=(x_4, \ldots, x_g, y_1,
\ldots, y_e)$. The maps $f$ and $\pi_{123}$ fit in a commutative
diagram, where $p_1: (\PP^1)^{g-3+e}\rightarrow (\PP^1)^{g-3}$ is
the projection on the first $g-3$ factors:
$$\begin{array}{ccc}

     \cc_{0, g}^e&    \stackrel{f}\longrightarrow  & (\PP^1)^{g-3+e}=\mm_{0, 4}^{g-3+e} \\
      \rmapdown{\pi_e} & \; & \rmapdown{p_1}  \\
    \mm_{0, g} &     \stackrel{\pi_{123}}\longrightarrow &  (\PP^1)^{g-3}=\mm_{0, 4}^{g-3} \\
\end{array}$$
Finally, for $2\leq k\leq e$ we define the diagonal loci
$\Delta_k\subset (\PP^1)^{g-3+e}$ as consisting of those points
$(x_4, \ldots, x_g, y_1, \ldots, y_e)$ for which at least $k$ of the
points $y_1, \ldots, y_e$ coincide. We need the following result concerning
existence of sublinear limit linear series of a fixed limit
$\mathfrak g^r_d$, having prescribed vanishing sequence at a given
point:
\begin{lemma}\label{inductive}
Let $X$ be a curve of compact type, $Y\subset X$ an irreducible
component and let $p\in Y$ be a smooth point of $X$. Assume that $l$
is a (refined) limit $\mathfrak g^r_d$ on $X$ and let $(a_0< a_1<
\ldots < a_r)$ be the vanishing sequence $a^l(p)$. We fix a
subsequence $(a_{j_0}<a_{j_1}<\ldots <a_{j_b})$ of $a^l(p)$, where
$0\leq b\leq r$.  Then there exists a limit $\mathfrak g_d^b$ on
$X$, say $ l'\subset l$, such that $a^{l'}(p)=(a_{j_0}, \ldots,
a_{j_b})$.
\end{lemma}
\begin{proof}
 Let us denote by $l:=\{l_Z=(L_Z, V_Z)\}_{Z\subset X}$ the original
 limit $\mathfrak g^r_d$ on $X$. For each integer $0\leq k\leq b$ there
 exists a section $\sigma_{j_k}\in V_Y$ such
 that $\mbox{ord}_p(\sigma_{j_k})=a_{j_k}$. We consider the subspace
 $W_Y:=<\sigma_{j_0}, \ldots, \sigma_{j_b}>\subset V_Y$. Since $\#\{\mbox{ord}_p(\sigma)\}_{\sigma \in W_Y}=b+1$,
 we obtain that $\mbox{dim}(W_Y)=b+1$ and we set $l'_Y:=(L_Y, W_Y)\in G^b_d(Y)$. Suppose now that $Z$ is a
  component of $X$ meeting $Y$ in a point $q$. We denote by $(c_{j_0}<c_{j_1}<\ldots <c_{j_b})$ the vanishing
  sequence $a^{l'_Y}(q)$. Let $(e_{j_0}<e_{j_1}<\ldots <e_{j_b})$ be the complementary sequence, that is,
   $e_{j_k}=d-c_{j_{b-k}}$ for each $0\leq k\leq b$. Then we can choose a section $\tau_k\in V_Z$ such that
 $\mbox{ord}_q(\tau_k)=e_{j_k}$. We define $W_Z:=<\tau_0, \ldots, \tau_b>\subset V_Z$. Because all
the entries $(e_{j_k})_{k=0}^b$ are distinct, we get that
$\mbox{dim}(W_Z)=b+1$ and then set $l'_Z:=(L_Z, W_Z)\in G^b_d(Z)$.
We continue inductively, and for each irreducible component
$Z'\subset X$ we obtain an aspect $l_{Z'}'=(L_{Z'}, W_{Z'})\in
G^b_d(Z')$. The collection $\{l_Z'\}_{Z\subset X}$ is the desired
limit $\mathfrak g_d^b$ on $X$.
\end{proof}
Next we explain how the assumption that for every $[C]\in \cM_g$
there exists a linear series $l\in G^r_d(C)$ with $V^{e-f}_e(l)\neq
\emptyset$, can be used to construct a flag curve $\tilde{R}\in
j(\mm_{0, g})$ such that all the $e$ points coming from the limit of
an effective divisor $D\in V_e^{e-f}(l)$ specialize to a connected
subcurve of $\tilde{R}$ having arithmetic genus at most
$\mbox{min}\{g, e\}$.
\begin{proposition}\label{schub}
Let $U\subset \cc_{0, g}^e$ be an irreducible component of the
closure of the locus of limits of $e$-secant divisors with respect
to linear series $\mathfrak g^r_d$ on flag curves from $\mm_g$.
Assuming that $\mathrm{dim}(U)=g-3+e-m$ \  with $0\leq m\leq e$,
there exists a point $([R, x_1, \ldots, x_g], \tilde{y}_1, \ldots,
\tilde{y}_e) \in W\cap(\mm_{0, g}\times_{\mm_g}
\overline{\mathcal{C}}_g^e)$ corresponding to a genus $g$ flag curve
$$\tilde{R}=R\cup_{x_1} E_1\cup \ldots \cup_{x_g} E_g \ \mbox{ and points } \tilde{y}_1, \ldots, \tilde{y}_e\in \tilde{R},$$ such that
either (i) \ $\tilde{y}_1=\cdots =\tilde{y}_e\in R-\{x_1, \ldots,
x_g\}$, or else, (ii)\  all the points $\tilde{y}_1, \ldots,
\tilde{y}_e$ lie on a connected subcurve $Y\subset \tilde{R}$
satisfying $p_a(Y)\leq \mathrm{min}\{m, g\}$ and $\#(Y\cap
(\overline{\tilde{R}-Y}))\leq 1$.
\end{proposition}
\begin{proof}
We start by noting that if $m=0$ then $U=\cc_{0, g}^e$ and
possibility (i) is satisfied. Thus we may assume that $m\geq 1$.
First, we claim that $\mbox{dim } f(U)=\mbox{dim }U=g-3+e-m$.
Indeed, since $\pi_e(U)=\mm_{0, g}$ it follows that
$p_1(f(U))=(\PP^1)^{g-3}$ and we choose a general point $t=(x_4,
\ldots, x_g)\in (\PP^1-\{0, 1, \infty\})^{g-3}$, such that $x_i\neq
x_j$ for $i\neq j$. Then $\pi_e^{-1}(t)=(\PP^1)^e$ and $f_{|
\pi_e^{-1}(t)}$ is an isomorphism onto its image, hence $f_{| U}$ is
birational onto its image as well. Obviously, when $m\geq g$ we can
take $Y=\tilde{R}$. From now on we shall assume that $1\leq m\leq
g-1$.

Let us assume first that $f(U)\cap \Delta_e\neq \emptyset$. Then
$\mbox{dim}\bigl( f(U)\cap \Delta_e\bigr)\geq g-m-2$. For dimension
reasons, there must exist a point $z=(x_4, \ldots, x_g, y_1, \ldots,
y_1)\in f(U)\cap \Delta_e$ such that either (i) at least $g-m-3$ of
the points $x_j$ with $4\leq j\leq g$ are mutually distinct and
belong to the set $\PP^1-\{0, 1, \infty, y_1\}$ and $y_1\in
\PP^1-\{0, 1, \infty\}$, or (ii) at least $g-m-2$ of the $x_j$'s
($4\leq j\leq g$) are mutually distinct and belong to the set
$\PP^1-\{0, 1, \infty, y_1\}$ and then $y_1\in \PP^1$ may, or may
not be equal to one of the points $0, 1$ or $\infty$. Suppose we are
in situation (i), the remaining case being similar.
 We fix a point $([R, x_1, \ldots,
x_g], y_1, \ldots, y_e)\in f^{-1}(z)$, hence $y_1, \ldots, y_e\in
R$. If $Z\subset R$ denotes the minimal connected subcurve of $R$
containing  all the points $y_1, \ldots, y_e$, then   $x_1, x_2,
x_3\in R-Z$, unless $y_1=\cdots =y_e$. (In the latter case either
$y_1\in R-\{x_1, \ldots, x_g\}$ which corresponds to the situation
when all the points $\tilde{y}_i=y_i$ specialize to the same smooth
point of $\tilde{R}$ lying on the rational spine, or else, if
$y_1=x_j$ for some $4\leq j\leq g$, then we can find a connected
subcurve of $\tilde{R}$ of genus $1$ containing $\tilde{y}_1,
\ldots, \tilde{y}_e$, where $p_R(\tilde{y}_i)=y_i$ for $1\leq i\leq
e$). Since at least $g-m=3+(g-m-3)$ of the points $x_1, \ldots, x_g$
lie on $Z^c$, it follows that  $\tilde{y}_1, \ldots, \tilde{y}_e$
lie on a connected subcurve of $\tilde{R}$ of genus $\leq m$, which
completes the proof in this case.

We are left with the possibility $f(U)\cap
\Delta_e=\emptyset$ and we denote by $k\leq e-1$ the largest integer
for which $f(U)\cap \Delta_k\neq \emptyset$ and by $L$ an
irreducible component of $f(U)\cap \Delta_{k}$. Since by definition
$f(U)\cap \Delta_{k+1}=\emptyset$, it follows that there exists a
point $t_0=(p_1, \ldots, p_e)\in (\PP^1)^{e}$ such that $L\subset
(\PP^1)^{g-3}\times \{t_0\}$. In particular, the projection map
$p_{1 | L}:L\rightarrow p_1(L)$ is $1:1$ and then $\mbox{dim }
p_1(L)=\mbox{dim}(L)\geq g-m+(e-k-2)\geq g-m$, unless $k=e-1$, when
$\mbox{dim }p_1(L)\geq g-m-1$. In the first case it follows that
there exists a point $(x_4, \ldots, x_g, p_1, \ldots, p_e)\in
f(U)\cap \Delta_k$ such that at least $g-m$ of the points $x_4,
\ldots, x_g$ are equal to a fixed point $r\in \PP^1-\{p_1, \ldots,
p_e\}$. In the second case, that is, when $k=e-1$, since
$\#\{p_i\}_{i=1}^e=2$, one of the points $0, 1$ or $ \infty$, say
$0$, does not appear among the $p_i$'s. Then we can find a point
$(x_4, \ldots, x_g, p_1, \ldots, p_e)\in f(U)\cap \Delta_{e-1}$ with
at least $g-m$ of the $x_j$'s equal to $0$.

The conclusion in both cases is that there exists a point $\bigl([R,
x_1, \ldots, x_g], y_1, \ldots, y_e\bigr)\in W\cap(\mm_{0,
g}\times_{\mm_g} \cc_g^e)$ corresponding to the flag curve
$\tilde{R}=R\cup_{x_1} E_1 \cup \ldots \cup_{x_g} E_g$, such that
the points $y_1,\ldots, y_e$ lie on a connected subcurve $Y\subset
\tilde{R}$ where $\#(Y\cap (\overline{\tilde{R}-Y}))\leq 1$ and
$p_a(Y)\leq m\leq e$.
\end{proof}

\noindent \emph{Proof of Theorem \ref{secant}.} We choose
$\tilde{R}=R\cup_{x_1} E_1\cup \ldots \cup_{x_g} E_g$ as above and
denote by $Y\subset \tilde{R}$ a connected subcurve onto which the
points $y_1, \ldots, y_e$ specialize. We know that either (a) \
$p_a(Y)=m\leq \mathrm{min}\{e, g\}$, or  (b) \ $y_1=\cdots =y_e\in
R-\{x_1, \ldots, x_g\}$.

We first deal with case (a) and dispose of (b) at the end using
\cite{EH2}. If $m<g$ we set $Z:=\overline{\tilde{R}-Y}$ and
$\{p\}:=Y\cap Z$ and we denote by $Y'$ and $Z'$ the components of
$Y$ and $Z$ respectively, containing the point $p$. When $m=g$, then
necessarily $e\geq g$ and $Y:= \tilde{R}, Z=\emptyset$ and $p\in
\tilde{R}$ is a general (smooth) point.  By assumption, $[\tilde{R},
y_1, \ldots, y_e]\in W$, hence there exists a proper flat morphism
$\phi:\mathcal{X}\rightarrow B$ satisfying the following properties:

\noindent $\bullet$  $\mathcal{X}$ is a smooth surface, $B$ is a
smooth affine curve, $0\in B$ is a point such that $\phi^{-1}(0)$ is
a curve stably equivalent to $\tilde{R}$ and $X_t=\phi^{-1}(t)$ is a
smooth projective curve of genus $g$ for $t\neq 0$. Moreover, there
are $e$ sections $\sigma_i:B\rightarrow \mathcal{X}$ of $\phi$
satisfying the condition $\sigma_i(0)=y_i\in \phi^{-1}(0)_{reg}$ for
all $1\leq i\leq e$.

\noindent $\bullet$ If $X_{\eta}:=\mathcal{X}-\phi^{-1}(0)$, then
there exists a line bundle $L_{\eta}\in \mbox{Pic}(X_{\eta})$ of
relative degree $d$ and a subvector bundle  $V_{\eta}\subset
\phi_*(L_{\eta})$ having rank $r+1$, such that for $t\neq 0$ we have
that $$\mbox{dim }V_t\cap H^0\bigl(X_t, L_t(-\sum_{j=1}^e
\sigma_j(t))\bigr)= r+1-e+f. $$ After possibly making a finite base
change and resolving the resulting singularities, the pair
$(L_{\eta}, V_{\eta})$ induces a (refined) limit $\mathfrak g^r_d$
on $\tilde{R}$, which we denote by $\mathfrak l$. The vector bundle
$V_{\eta}\cap \phi_*\bigl(L_{\eta}\otimes
\OO_{X_{\eta}}(-\sum_{j=1}^e \sigma_j(B-\{0\}))\bigr)$ induces a
limit linear series $\mathfrak g^{r-e+f}_{d-e}$ on $\phi^{-1}(0)$
which we denote by $\mathfrak m$. For a component $A$ of
$\phi^{-1}(0)$, if $(L_A, V_A)\in G^r_d(A)$ denotes the $A$-aspect
of $\mathfrak l$, then there exists a unique \emph{effective}
divisor $D_A\in A_e$ supported only at the points from $(A\cap
\bigcup_{j=1}^e \sigma_j(B))\bigcup (A\cap
\overline{\phi^{-1}(0)-A})$ such that the $A$-aspect of $\mathfrak
m$ is of the form
$$\mathfrak m_A=\bigl(M_A:=L_A\otimes \OO_A(-D_A), \ W_A\subset V_A\cap H^0(M_A)\bigr)\in
G^{r-e+f}_{d-e}(A).$$ The collection $\mathfrak m_Y:=\{\mathfrak
m_A\}_{A\subset Y}$ forms a limit $\mathfrak g^{r-e+f}_{d-e}$ on
$Y$. We denote by $(a_0< \ldots <a_r)$ the vanishing sequence of
$\mathfrak l_{Y'}$ at $p$, thus $\{a_i
\}_{i=0}^r=\{\mbox{ord}_p(\sigma)\}_{\sigma \in V_{Y'}}$ and we
denote by $(b_0<\ldots <b_r)$ the vanishing sequence $a^{\mathfrak
l_{Z'}}(p)$. By ordering the set $\{\mbox{ord}_p(\sigma)\}_{\sigma
\in W_{Y'}}$ we obtain a subsequence $(a_{i_0}<\ldots
<a_{i_{r-e+f}})$ of $a^{\mathfrak l_{Y'}}(p)$. When we order the
entries in $\{a_i\}_{i=0}^r-\{a_{i_k}\}_{k=0}^{r-e+f}$ we obtain a
new sequence $(a_{j_0}<a_{j_1}<\ldots <a_{j_{e-f-1}})$. Using Lemma
\ref{inductive}, we find that there exists a limit linear series
$\mathfrak l'_Y$ of type $\mathfrak g^{e-f-1}_{d}$ on $Y$ with the
property that $a^{\mathfrak l'_Y}(p)=(a_{j_0}, a_{j_1}, \ldots,
a_{j_{e-f-1}})$.

Let us assume first that we are in the situation $m<g$, hence $Z\neq
\emptyset$. The point $p\in Y$ lies on a rational component which
implies the following inequality corresponding to $Y$ (see also
\cite{EH2}, Theorem 1.1):
\begin{equation}\label{egy1}
V_1:=\rho(m, e-f-1, d)-\sum_{k=0}^{e-f-1} a_{j_k}+{e-f\choose 2}\geq
0.
\end{equation}
Applying the same principle for the limit linear series $\mathfrak
m_Y$ on $Y$, we find that the adjusted Brill-Noether number with
respect to the point $p$ is non-negative:
\begin{equation}\label{egy2} V_2:=\rho(m, r-e+f,
d-e)-\sum_{k=0}^{r-e+f} a_{i_k}+{r+1-e+f\choose 2}\geq 0.
\end{equation}
 Next we turn our attention to $Z$ and use the fact that
 the point $p\in Z$ does not lie on an elliptic component, hence
 $[Z, p]$ satisfies the "strong" pointed Brill-Noether theorem:
 \begin{equation}\label{egy3}
 V_3:=\rho(g-m, r, d)-\sum_{k=0}^r b_k+{r+1\choose 2}\geq 0.
 \end{equation}
 If we add (\ref{egy1}), (\ref{egy2}) and (\ref{egy3}) together and
 use that $\sum_{k=0}^r b_k+\sum_{k=0}^{r-e+f}
 a_{i_k}+\sum_{k=0}^{e-f-1} a_{j_k}=(r+1)d$, we obtain the
 inequality
 $$\rho(g, r, d)-f(r+1-e+f)+e\geq e-m\geq 0.$$

 The case $m=g$, when $Y=\tilde{R}$, is similar but simpler. We add together
 (\ref{egy1}) and (\ref{egy2}) (now there is no
 (\ref{egy3})) and we write the following inequalities:
 $$\rho(g, r, d)+e-f(r+1-e+f)=\Bigl(\rho(g, r-e+f,
 d-e)-\sum_{k=0}^{r-e+f} a_{i_k}+{r+1-e+f\choose
 2}\Bigr)+$$
 $$+\Bigl(\rho(g, e-f-1, d)-\sum_{k=0}^{e-f-1}a_{j_k}+{e-f\choose
 2}\Bigr)+\sum_{k=0}^{r-e+f} a_{i_k}+\sum_{k=0}^{e-f-1}
 a_{j_k}-{r+1\choose 2}+e-g\geq e-g\geq 0,$$
 since $\sum_{k=0}^{r-e+f} a_{i_k}+\sum_{k=0}^{e-f-1} a_{j_k}\geq {r+1\choose 2}$. Thus  we
 obtain the same numerical conclusion as in the case $m<g$.

 Assume now that we are in the case (b) when $y_1=\cdots =y_e\in
 R-\{x_1, \ldots, x_g\}$. Then reasoning as above, we find a limit
 $\mathfrak g^r_d$ on $\tilde{R}$ having vanishing sequence at $y_1$
 at least $(0, 1, \ldots, e-f-1, e, e+1, \ldots, r+f-1, r+f)$. Using once
 more \cite{EH2}, Theorem 1.1, we obtain the inequality
$$\rho(g, r, d)+e-f(r+1-e+f)\geq \rho(g, r, d)-f(r+1-e+f)\geq 0.$$

 Using the semicontinuity of the dimension of the fibres,
 it follows that for a general curve $[C]\in \cM_g$, if $\pi_1:\Sigma_C\rightarrow C_e$ is the first
 projection, then the minimal fibre dimension of $\pi_1$ cannot exceed
 the dimension of the space of pairs of limit linear series
 $\mathfrak l\supset  \mathfrak m$
 consisting of a $\mathfrak g^r_d\supset \mathfrak g_{d-e}^{r-e+f}$ on the flag curve
 $\phi^{-1}(0)$ such that $\mathfrak m=\mathfrak l(-D_e)$, where $D_e$ is a degree $e$ effective divisor
 on $\phi^{-1}(0)$ with the property that $\mbox{supp}(D_e)\subset Y\cap \phi^{-1}(0)_{reg} $. Since the map
 $(\mathfrak l\supset \mathfrak m, \mathfrak m_Y, \mathfrak l_Y')\mapsto (\mathfrak m_Y, \mathfrak l_Y', \mathfrak l_Z)\in \tilde{G}_{d-e}^{r-e+f}(Y)\times \tilde{G}_{d}^{e-f-1}(Y)\times \tilde{G}_{d}^r(Z)$ is injective, it follows that  for a general
 divisor $D_{gen}\in \pi_1(\Sigma_C)$ we have the estimate
 $$\mbox{dim }\pi_1^{-1}(D_{gen})\leq
 V_1+V_2+V_3=\rho(g, r, d)-f(r+1-e+f)+m,$$ hence
 $\mbox{dim}(\Sigma_C)=\mbox{dim }\pi_1^{-1}(D_{gen})+e-m\leq \rho(g,
 r,d)-f(r+1-e+f)+e$. This finishes the proof of Theorem \ref{secant}.
$\hfill$ $\Box$

\section{Existence of linear series with secant planes}

We turn our attention to showing existence of linear series which
possess $e$-secant $(e-f-1)$-planes. The strategy we pursue is to
construct limit linear series $\mathfrak g^r_d$ on a curve of
compact type $[Y\cup_p Z]\in \mm_g$, where $(Y, p)$ and $(Z, p)$ are
suitably general smooth pointed curves of genus $e$ and $g-e$
respectively. These limit $\mathfrak g^r_d$'s will carry  a
sublinear series $\mathfrak g_{d-e}^{r-e+f}=\mathfrak g^r_d(-D_e)$,
where $D_e$ is a degree $e$ effective divisor on $Y$. Like in the
proof of Theorem \ref{secant}, such $\mathfrak g^r_d$'s are
determined by their $Z$-aspect and by a pair of linear series
$(\mathfrak g_{d-e}^{r-e+f}, \mathfrak g_d^{e-f-1})$ on $Y$. We
determine the dimension of the space of such pairs, which will
enable us to show that the original pair $(\mathfrak
g_{d-e}^{r-e+f}, \mathfrak g_d^{e-f-1})$ on $Y\cup_p Z$ can be
smoothed to every smooth curve of genus $g$. This will finish the
proof of Theorem \ref{existence}.

We start by choosing two general pointed curves $[Y, p]\in \cM_{e,
1}$ and $[Z, p]\in \cM_{g-e, 1}$ such that both $(Y, p)$ and $(Z,
p)$ satisfy the Brill-Noether theorem with prescribed ramification
(cf. \cite{EH2}, Theorem 1.1 and Proposition 1.2): If
$\overline{\alpha}: 0\leq \alpha_0\leq \ldots \leq \alpha_r\leq d-r$
is a Schubert index of type $(r, d)$, then $(Y, p)$ possesses a
$\mathfrak g^r_d$ with ramification sequence $\geq
\overline{\alpha}$ at the point $p$, if and only if
\begin{equation}\label{eisharris}
\sum_{i=0}^r \mathrm{max}\{\alpha_i+g(Y)-d+r, 0\}\leq g(Y).
\end{equation}
In case this inequality is satisfied, then $\mbox{dim } G^r_d(Y, p,
\overline{\alpha})=\rho(g, r, d, \overline{\alpha})$ (One obviously
has a similar statement for $[Z, p]$).

We denote by $\pi:\mathcal{X}\rightarrow (T, 0)$ the versal
deformation space of the stable curve $\pi^{-1}(0)=X_0:=Y\cup_p Z$.
Let $\Delta\subset T$ be the boundary divisor corresponding to
singular curves, and we write
$\pi^{-1}(\Delta)=\Delta_e+\Delta_{g-e}$, where $\Delta_e$ (resp.
$\Delta_{g-e}$) is the divisor corresponding to the marked point
lying on the component of genus $e$ (resp. $g-e$). We consider the
$e$-fold fibre product
$\mathcal{U}:=(\mathcal{X}-\Delta_{g-e})\times_T \cdots \times_T
(\mathcal{X}-\Delta_{g-e})$, the projection
$\phi:\mathcal{U}\rightarrow T$ and the induced curve
$p_2:\mathcal{X}\times_T \mathcal{U}\rightarrow \mathcal{U}$. Then
we introduce the stack of limit linear series of type $\mathfrak
g^r_d$ over $\mathcal{U}$
$$\sigma:\widetilde{\mathfrak{G}}_d^r(\mathcal{X}\times_T \mathcal{U}/\mathcal{U})\rightarrow
\mathcal{U}, \mbox{ where
}\widetilde{\mathfrak{G}}_d^r(p_2)=\widetilde{\mathfrak{G}}_d^r(\mathcal{X}\times_T
\mathcal{U}/\mathcal{U})=\widetilde{\mathfrak{G}}^r_d(\pi)\times_T
\mathcal{U},$$ and we write $\tau:=\phi\circ
\sigma:\widetilde{\mathfrak{G}}_d^r(p_2)\rightarrow T$ (see
\cite{EH1} Theorem 3.4, for details on the construction of
$\widetilde{\mathfrak G}^r_d(\pi)$). The fibre $\tau^{-1}(t)$
corresponding to a point $t\in \Delta$ (in which case one can write
 $\pi^{-1}(t)=Y_t\cup Z_t$, with $g(Y_t)=e, g(Z_t)=g-e$),
parameterizes limit $\mathfrak g^r_d$'s on $Y_t\cup Z_t$ together
with $e$-tuples $(x_1, \ldots, x_e)\in (Y_t-Y_t\cap Z_t)^e$. Let us
denote by $\mathcal{L}_Y$ a degree $d$ Poincar\'e bundle on
$\pi_2:\mathcal{X}\times _T
\widetilde{\mathfrak{G}}_d^r(p_2)\rightarrow \widetilde{\mathfrak
G}_d^r(p_2)$ characterized by the property that its restriction to
curves of type $Y_t\cup Z_t$ are line bundles of  bidegree $(d, 0)$.
We also write $\mathcal{V}_Y\subset (\pi_{2})_* (\mathcal{L}_Y)$ for
the rank $r+1$ tautological bundle whose fibres correspond to the
global sections of the genus $e$-aspect of each limit $\mathfrak
g^r_d$. Finally, for $1\leq j\leq e$, we denote by $D_j\subset
\mathcal{X}\times_T \widetilde{\mathfrak G}^r_d(p_2)$ the diagonal
divisor corresponding to pulling back the diagonal under the map
$\mathcal{X}\times_T \widetilde{\mathfrak G}^r_d(p_2)\rightarrow
\mathcal{X}\times_T \mathcal{X}$ which projects onto the $j$-th
factor, that is, $(x, l, x_1, \ldots, x_e)\mapsto (x, x_j)$ where
$x, x_1, \ldots, x_e\in \pi^{-1}(t)$. There exists an evaluation
vector bundle morphism over $\widetilde{\mathfrak G}^r_d(p_2)$
$$\chi:\mathcal{V}_Y\rightarrow (\pi_{2})_*(\mathcal{L}_Y\otimes
\OO_{\sum_{j=1}^e D_j})$$ and we denote by $\H$ the rank $e-f$
degeneracy locus of the map $\chi$. Set-theoretically, $\H$ consists
of those points $(t, l, x_1, \ldots, x_e)$ with $\phi(x_1, \ldots,
x_e)=t\in T$ and $l\in \widetilde{G}^r_d(\pi^{-1}(t))$, satisfying
the condition that $\mbox{dim } l(-x_1-\cdots -x_e)\geq r+1-e+f$.
The dimension of every irreducible component of $\H$ is at least $
\rho(g, r, d)+\mbox{dim }T+e-f(r+1-e+f)$.

 In order to show that $\tau:\H\rightarrow T$ is
dominant, it suffices to prove that $\tau^{-1}(0)$ has at least one
irreducible component of dimension $\rho(g, r, d)+e-f(r+1-e+f)$.
This in fact will prove the stronger statement that $\Sigma_C\neq
\emptyset$ for \emph{every} $[C]\in \cM_g$. Indeed, even though
$\tau:\widetilde{\mathfrak G}^r_d(p_2)\rightarrow T$ is not a proper
morphism, the restriction
$\tau_{\tau^{-1}(T-\Delta)}:\tau^{-1}(T-\Delta)\rightarrow T-\Delta$
is proper, hence there exists an irreducible component of $\H$ which
maps onto $T-\Delta$. Since $\pi:\mathcal{X}\rightarrow (T, 0)$ can
be chosen in such a way that there exists a point $t\in T$ with
$\pi^{-1}(t)\cong C$, this proves our contention. We set the integer
$$\alpha_0:=\Bigl[\frac{\rho(e, r-e+f,
d-e)}{r+1-e+f}\Bigr]=\Bigl[\frac{e}{r+1-e+f}\Bigr]+d-r-f-e,$$ thus
we can write $\rho(e, r-e+f, d-e)=\alpha_0\cdot (r+1-e+f)+c$, where
$0\leq c\leq r-e+f$. Then there exists a unique Schubert index of
type $(r-e+f, d-e)$, $$\overline{\alpha}:0\leq \alpha_0\leq
\alpha_1\leq \ldots \leq \alpha_{r-e+f} \leq d-r-f,$$ with
$\alpha_{r-e+f}-\alpha_0\leq 1$, such that $\sum_{j=0}^{r-e+f}
\alpha_j=\rho(e, r-e+f, d-e)$. We have that $\alpha_j=\alpha_0$ for
$0\leq j\leq r-e+f-c$ and $\alpha_j=\alpha_0+1$ for $r-e+f-c+1\leq
j\leq r-e+f$. Note that since
$\alpha_0+g(Y)-(d-e)+r-e+f=[e/(r+1-e+f)]\geq 0$, condition
(\ref{eisharris}) is verified and the variety $G^{r-e+f}_{d-e}(Y, p,
\overline{\alpha})$ is non-empty of dimension $\rho(e, r-e+f,
d-e)-\sum_{j=0}^{r-e+f}\alpha_j=0$.

Next we set $\beta_0:=[e/(e-f)]$ and write \ $e=\beta_0\cdot
(e-f)+\tilde{c}$, where $0\leq \tilde{c} \leq e-f-1.$ Then there
exists a unique Schubert index of type $(e-f-1, 2e-f-1)$
$$\overline{\beta}:0\leq \beta_0\leq \beta_1\leq \ldots \leq
\beta_{e-f-1}\leq e,$$ such that $\beta_{e-f+1}-\beta_0\leq 1$ and
$\sum_{j=0}^{e-f-1} \beta_j=e$. Precisely, $\beta_j=\beta_0$ for
$0\leq j\leq e-f-\tilde{c}-1$ and $\beta_j=\beta_0+1$ for
$e-f-\tilde{c}\leq j\leq e-f-1$. By (\ref{eisharris}), the variety
$G^{e-f-1}_{2e-f-1}(Y, p, \overline{\beta})$ is non-empty and of
dimension $e-\sum_{j=0}^{e-f-1}\beta_j=0$.

First we are going to prove Theorem \ref{existence} under the
assumption that there exist two linear series $(A, W_A)\in
G^{r-e+f}_{d-e}(Y, p, \overline{\alpha})$ and $(L, W_L)\in
G^{e-f-1}_{2e-f-1}(Y, p, \overline{\beta})$ satisfying the condition
\begin{equation}\label{assumption2}
 H^0\bigl(Y, L\otimes A^{\vee}\otimes \OO_Y((d+f-2e)\cdot
 p)\bigr)=0.
\end{equation}
Note that $\mbox{deg}\bigl(L\otimes A^{\vee}\otimes
\OO_Y((d+f-2e)\cdot p\bigr)=g(Y)-1$, and (\ref{assumption2}) states
that a suitable translate of at least one of the finitely many line bundles
of type $L\otimes A^{\vee}$ lies outside the theta divisor of $Y$.

\begin{remark}
Condition (\ref{assumption2}) is a subtle statement concerning $[Y,
p]$. It is not true that (\ref{assumption2}) holds for \emph{every}
choice of $(A, W_A)\in G_{d-e}^{r-e+f}(Y, p, \overline{\alpha})$ and
$(L, W_L)\in G_{2e-f-1}^{e-f-1}(Y, p, \overline{\beta})$. For
instance, in the case $e=2r-2$ and $f=r-1$, corresponding to
$(2r-2)$-secant $(r-2)$-planes which every curve $Y\subset \PP^r$ is
expected to possess in finite number, we obtain that $A=B\otimes
\OO_Y((d-3r+2)\cdot p)$, where $B\in W^1_r(Y)$ and $L\otimes
\OO_Y(-2p)\in W^{r-2}_{3r-6}(Y)$. By Riemann-Roch, we can write that
$L=K_Y\otimes \OO_Y(2\cdot p)\otimes \tilde{B}^{\vee}$, where
$\tilde{B}\in W^1_r(Y)$ and then (\ref{assumption2}) translates into
the vanishing statement $H^0(Y, B\otimes \tilde{B}\otimes
\OO_Y(-3\cdot p))=0$. The curve $Y$ has $\frac{(2r-2)!}{r! (r-1)!}$
pencils $\mathfrak g^1_r$. If we choose $B\neq \tilde{B}\in
W^1_r(Y)$, then $h^0(Y, B\otimes \tilde{B})\geq 4$ and
(\ref{assumption2}) has no chance of being satisfied. If
$B=\tilde{B}$, then the Gieseker-Petri theorem implies that the map
$H^0(Y, B)\otimes H^0(Y, K_Y\otimes B^{\vee})\rightarrow H^0(Y,
K_Y)$ is an isomorphism, whence $h^0(Y, B^{\otimes 2})=3$. Choosing
$p\in Y$ outside the set of ramification points of the finitely many
line bundles $B^{\otimes 2}$ where $B\in W^1_r(Y)$, we obtain that
$H^0(B^{\otimes 2}\otimes \OO_Y(-3\cdot p))=0$. Therefore in this
case, condition $(\ref{assumption2})$ is equivalent to the
Gieseker-Petri theorem.
\end{remark}

We shall study when (\ref{assumption2}) is actually satisfied. We
note that by the Riemann-Roch theorem, (\ref{assumption2}) also
implies that $h^0\bigl(Y, L\otimes A^{\vee}\otimes
\OO_Y((d+f-2e+1)\cdot p)\bigr) =1$.
 Assuming that $(A, W_A)\in G^{r-e+f}_{d-e}(Y, p,
\overline{\alpha})$ and $(L, W_L)\in G^{e-f-1}_{2e-f-1}(Y, p,
\overline{\beta})$ satisfy (\ref{assumption2}), it follows from
Riemann-Roch that there exists a unique effective divisor of degree
$e$
$$D\in |L\otimes \OO_Y((d-2e+f+1)\cdot p)\otimes A^{\vee}|,$$ and
moreover $p\notin \mbox{supp}(D)$. We introduce the space of
sections
$$V_Y:=W_A+W_L\subset H^0\bigl(Y, L\otimes \OO_Y((d-2e+f+1)\cdot p)\bigr),\ \mbox{ where we view }$$
$$ W_A\subset H^0\bigl(L\otimes \OO_Y((d-2e+f+1)\cdot p-D)\bigr) \mbox{ and } W_L\subset H^0(L)\subset
H^0\bigl(L\otimes \OO_Y((d-2e+f+1)\cdot p)\bigr).$$

We claim that $\mbox{dim}(V_Y)=r+1$, hence $\mathfrak{l}_Y=(L\otimes
\OO_Y((d-2e+f+1)\cdot p), V_Y)\in G^r_d(Y)$. Moreover,
$\mathfrak{l}_Y$ has the following vanishing sequence at $p$:
\begin{equation}\label{vanseq1}
a^{\mathfrak{l}_Y}(p)=(\alpha_0, \ldots, \alpha_{r-e+f}+r-e+f, \beta_0+d-2e+f+1,\beta_1+d-2e+f+2, \ldots,
\beta_{e-f-1}+d-e).
\end{equation}
Indeed, our original assumption  $f(r+1-e+f)\geq e$ is equivalent
with the inequality $\alpha_{r-e+f}+r-e+f<d-2e+f+1$, which shows
that the sequence (\ref{vanseq1}) contains $r+1$ distinct entries.
Since $p\notin \mbox{supp}(D)$, we obtain that the vanishing orders
of the sections from $W_A\subset H^0(L\otimes \OO_Y((d-2e+f+1)\cdot
p))$ are precisely
$$\alpha_0, \alpha_1+1, \ldots, \alpha_{r-e+f}+r-e+f,$$ while those of
the sections from $W_L\subset H^0(L\otimes \OO_Y((d-2e+f+1)\cdot
p))$ are precisely $$\beta_0+d-2e+f+1, \beta_1+d-2e+f+2, \ldots,
\beta_{e-f-1}+e-f-1+d-2e+f+1=\beta_{e-f-1}+d-e.$$ We have found
$r+1$ sections from $V_Y$ having distinct vanishing orders at the
point $p$, hence $\mbox{dim}(V_Y)=r+1$. Moreover,  $a^{{\mathfrak
l}_Y}(p)$ is equal to the sequence (\ref{vanseq1}).

 Next we choose a
linear series $\mathfrak l_Z\in G^r_d(Z, p)$ such that $\{\mathfrak
l_Y, \mathfrak l_Z\}$ is a refined limit $\mathfrak g^r_d$. Then the
ramification sequence of $\mathfrak l_Z$ at the point $p$ must be
equal to
$$\alpha^{\mathfrak l_Z}(p)=\overline{\gamma}:=(e-\beta_{e-f-1}, e-\beta_{e-f-2},  \ldots, e-\beta_0, d-r-\alpha_{r-e+f},
\ldots, d-r-\alpha_1, d-r-\alpha_0).$$ We claim that condition
(\ref{eisharris}) is satisfied for $Z$ and that the variety $G^r_d(Z, p,
\overline{\gamma})$ is non-empty and of dimension $\rho(g-e, r, d,
\overline{\gamma})=\rho(g, r, d)+e-f(r+1-e+f)$. For this to happen,
one has to check that the following inequality holds:
\begin{equation}\label{nou} \sum_{j=0}^r
\mathrm{max}\{\alpha^{\mathfrak l_Z}_j(p)+g-e-d+r, 0\}\leq g-e.
\end{equation}
 There are two things to notice: First, that by direct
computation we have that
 $$\alpha^{\mathfrak
l_Z}_{e-f}(p)+g-e-d+r=g-e-\alpha_{r-e+f}=(g-d+r)+\bigl[f-\frac{e}{r+1-e+f}\bigr]\geq
0,$$ hence $\alpha^{\mathfrak l_Z}_{j}(p)+g-e-d+r\geq 0$ for all
$e-f\leq j\leq r$. Second, that since $0\leq
\beta_{e-f-1}-\beta_0\leq 1$, in order to estimate the sum of the
first $e-f$ terms in the sum (\ref{nou}), there are two cases to consider.
Either $\alpha^{\mathfrak l_Z}_0(p)+g-e-d+r\geq 0$, in which case we
find that $$\sum_{j=0}^r \mathrm{max}\{\alpha^{\mathfrak
l_Z}_j(p)+g-e-d+r, 0\}=\sum_{j=0}^r (\alpha^{\mathfrak
l_Z}_j(p)+g-e-d+r)=$$ $$=g-e-\rho(g-e, r, d, \overline{\gamma})=
g-e-\bigl(\rho(g, r, d)+e-f(r+1-e+f)\bigr)\leq g-e.$$ Else, if
$\alpha^{\mathfrak l_Z}_0(p)+g-e-d+r\leq -1$, then also
$\alpha^{\mathfrak l_Z}_j(p)+g-e-d+r\leq 0$ for $0\leq j\leq e-f-1$
and the left hand side of (\ref{nou}) equals
$$\sum_{j=e-f}^{r}(\alpha_j^{\mathfrak l_Z}(p)+g-e-d+r)=
(r+1-e+f)(g-e)-\sum_{i=0}^{r-e+f} \alpha_{i}=g-e-\rho(g, r-e+f, d-e)\leq g-e.$$ In both cases
the inequality (\ref{eisharris}) is satisfied which proves our
claim.

Since the chosen $(A, W_A)\in G^{r-e+f}_{d-e}(Y, p,
\overline{\alpha})$ and $(L, W_L)\in G^{e-f-1}_{2e-f-1}(Y, p,
\overline{\beta})$ are isolated points in their corresponding
varieties of linear series on $Y$, it follows that limit $\mathfrak
g^r_d$'s on $X_0$ constructed in the way we just described, fill-up
a component of $\tau^{-1}(0)\subset \H$.

Indeed, suppose $(\mathfrak n_Y, \mathfrak n_Z, \tilde{D})\in \H$ is
a point lying in the same irreducible component of $\tau^{-1}(0)$ as
$(\mathfrak l_Y, \mathfrak l_Z, D)$. Here, $\mathfrak n_Y\in
G^r_d(Y),\ \mathfrak n_Z\in G^r_d(Z, p, \overline{\gamma})$ and
$\tilde{D}\in Y_e$ is a divisor such that $p\notin
\mbox{supp}(\tilde{D})$. Then $a^{\mathfrak n_Y}(p)=a^{\mathfrak
l_Y}(p)$ which is given by (\ref{vanseq1}), therefore $\mathfrak
n_Y(-(d-2e+f+1)\cdot p)\in G^{2e-f-1}_{e-f-1}(Y, p,
\overline{\beta})$ which is a reduced $0$-dimensional variety. This
implies that $\mathfrak n_Y(-(d-2e+f+1)\cdot p)=(L, W_L)$. Next, we
consider the linear series $\mathfrak n_Y(-\tilde{D})\in
G^{r-e+f}_{d-e}(Y)$. Since $p\notin \mbox{supp}(\tilde{D})$, the
vanishing sequence of this linear series is a subsequence of length
$r+1-e+f$ of $a^{\mathfrak l_Y}(p)$. Necessarily, $\alpha^{\mathfrak
n_Y(-\tilde{D})}(p)\geq \overline{\alpha}$ and because $\rho(e,
r-e+f, d-e, \overline{\alpha})=0$, we must have that $\mathfrak
n_Y(-\tilde{D})\in G^{r-e+f}_{d-e}(Y, p, \overline{\alpha})$ which
is a discrete set, hence $\mathfrak n_Y(-\tilde{D})=(A, W_A)$ and
$\tilde{D}=D\in Y_e$. This shows that $\mathfrak n_Y=\mathfrak l_Y$
and every point of this component of $\tau^{-1}(0)$ is determined by
the $\mathfrak n_Z$. The dimension of this component is thus equal
to
$$\rho(e, r-e+f, d-e, \overline{\alpha})+\rho(g-e, r, d,
\overline{\gamma})+\rho(e, e-f-1, 2e-f-1, \overline{\beta})=\rho(g,
r, d)-f(r+1-e+f)+e,$$ which finishes the proof of Theorem
\ref{existence}, subject to proving assumption (\ref{assumption2}).

\begin{remark} A slight variation of the argument described above, enables us to prove
Theorem \ref{existence} even in some cases when we cannot establish
(\ref{assumption2}). We start with a linear series $(A, W_A)\in
G^{r-e+f}_{d-e}(Y, p, \overline{\alpha})$ and assume that the
following condition holds:
\begin{equation}\label{assumption3}
H^0\bigl(Y, \OO_Y((d-1)\cdot p)\otimes A^{\vee}\bigr)=0.
\end{equation}
There exists a unique divisor $D\in |\OO_Y(d\cdot p)\otimes
A^{\vee})|$ and (\ref{assumption3}) guarantees that $p\notin
\mbox{supp}(D)$. We define the space of sections
$$V_Y:=H^0(\OO_Y(2e-f-1)\cdot p)+ W_A\subset H^0(\OO_Y(d\cdot p)), \ \mbox{ where
 } \ W_A\subset H^0( \OO_Y(d\cdot p-D)).$$
Reasoning along the same lines as in the previous case, since
$p\notin \mbox{supp}(D)$ we find that $\mbox{dim}(V_Y)=r+1$, hence
$\mathfrak l_Y=(\OO_Y(d\cdot p),
 V_Y)\in G^r_d(Y)$. Moreover, we can check that
$$
a^{\mathfrak l_Y}(p)=(\alpha_0, \alpha_1+1, \ldots,
\alpha_{r-e+f}+r-e+f, d-2e+f+1, d-2e+f+2, \ldots, d-e-1, d).
$$
Like in the previous situation, we choose a linear series $\mathfrak
l_Z\in G^r_d(Z, p)$ such that $\{\mathfrak l_Y, \mathfrak l_Z\}$ is
a refined limit $\mathfrak g^r_d$. Thus we must have the following
ramification sequence at $p$:
$$\alpha^{\mathfrak l_Z}(p)=\overline{\gamma}:=(0, e,  \ldots, e, d-r-\alpha_{r-e+f},
\ldots, d-r-\alpha_1, d-r-\alpha_0).$$  Condition (\ref{eisharris})
which guarantees the existence of $\mathfrak l_Z$ is satisfied if
and only if $$\rho(g, r,d)\geq f(r+1-e+f)-(g-d+r), \ \mbox{ in the
case }g-d+r<e
$$
and $$ \rho(g, r, d)\geq f(r+1-e+f)-e, \ \mbox{ in the case
}g-d+r\geq e. $$ Since we are always working under the hypothesis
$\rho(g, r, d)-f(r+1-e+f)+e\geq 0$, we see that the previous condition holds
whenever $g-d+r\geq e$, and that, in general, $\mathfrak l_Z\in
G^r_d(Z, p, \overline{\gamma})$ exists if and only if
\begin{equation}\label{ass4} \rho(g, r,d)\geq f(r+1-e+f)-(g-d+r).
\end{equation} Assuming (\ref{ass4}), the variety $G^r_d(Z, p,
\overline{\gamma})$ is non-empty of dimension $\rho(g-e, d, r,
\overline{\gamma})=\rho(g, r, d)-f(r+1-e+f)+e$. The same argument as
before shows that limit $\mathfrak g^r_d$'s on $X_0$ constructed in
such a way, fill-up a component of $\tau^{-1}(0)\subset \H$ of
expected dimension $\rho(g, r, d)-f(r+1-e+f)+e$, which finishes the
proof.
\end{remark}

Now we complete the proof of Theorem \ref{existence} by discussing
under which assumptions we can establish (\ref{assumption2}):

\noindent{\emph{Proof of Theorem \ref{existence}}}. We retain the
notation introduced above and show that there exist two linear
series $(A, W_A)\in G_{d-e}^{r-e+f}(Y, p, \overline{\alpha})$ and
$(L, W_L)\in G_{2e-f-1}^{e-f-1}(Y, p, \overline{\beta})$ satisfying
(\ref{assumption2}) whenever one of the following conditions is
satisfied: $$(i) \mbox{ } \  2f\leq e-1,\ \ (ii)\mbox{ } \ e=2r-2
\mbox{ and } f=r-1, \ \mbox{ } (iii)\mbox{ } \
 e<2(r+1-e+f).$$ As we already explained,  (\ref{assumption2}) in case $(ii)$ is a
 consequence of the Gieseker-Petri theorem.

We now treat case $(i)$ when $\beta_0=1$ and $\tilde{c}=f\leq
e-f-1$. By Riemann-Roch we find that $L=K_Y\otimes
\OO_Y((e-2f+2)\cdot p)\otimes B^{\vee}$, where $B\in W^1_{e-f+1}(Y)$
is a pencil such that $h^0\bigl(Y,B\otimes \OO_Y(-(e-2f+1)\cdot
p)\bigr)\geq 1$ (There are finitely many such $B\in W^1_{e-f+1}(Y)$
for a generic choice of $[Y, p]\in \cM_{e, 1}$). Applying the
base-point-free pencil trick, (\ref{assumption2}) is equivalent to
the injectivity of the multiplication map $$\mu_{B, M}: H^0(Y,
B)\otimes H^0(Y, M)\rightarrow H^0(Y, B\otimes M),$$ where
$M:=K_Y\otimes A^{\vee}\otimes \OO_Y((d-f-e+2)\cdot p)\in
W_{2e-f}^{e-f}(Y)$ is a complete linear series with vanishing
sequence at $p$ equal to
\begin{equation}\label{vanseq3}
a^M(p)=(0, 1, \ldots, e-f-a-1, e-f-a+c, r-a+2, r-a+3, \ldots, r,
r+1).
\end{equation}
Here we have set $a:=[e/(r+1-e+f)]$, hence we can write $e=a\cdot
(r+1-e+f)+c$, where $0\leq c\leq r-e+f$. By assumption we have that
$e-2a>c$ and clearly $\rho(M, \alpha^M(p))=0$, that is, there are
finitely many $M\in W^{e-f}_{2e-f}(Y)$ satisfying (\ref{vanseq3}).

To prove that $\mu_{B, M}$ is injective, we degenerate $[Y, p]\in
\cM_{e, 1}$ to a particular stable curve:
 $[Y_0, p_0]:=[E_0\cup_{p_1} E_1\cup \ldots\cup E_{e-2a-1}\cup_{p_{e-2a}} T, p_0]$, where $E_0,
 \ldots, E_{e-2a-1}$ are elliptic curves, $[T=E_{e-2a}, p_{e-2a}]\in \cM_{2a, 1}$ is a
 Petri general smooth pointed curve and the points $p_i, p_{i+1}\in E_i$ are such that $p_{i+1}-p_i\in
\mbox{Pic}^0(E_i)$ is not a torsion class for $0\leq i\leq e-2a-1$.
Note that $p_0$ lies on the first component $E_0$. By contradiction,
we assume that $\mu_{B, M}$ is not injective for every $[Y, p]\in
\cM_{e, 1}$ and for each of the finitely many linear series $M\in
W_{2e-f}^{e-f}(Y)$ satisfying (\ref{vanseq3}) and each $B\in
G^1_{e-f+1}\bigl(Y, p, (0, e-2f)\bigr)$.

We construct a limit $\mathfrak g_{2e-f}^{e-f}$ on $[Y_0, p_0]$, say
$\mathfrak m=\{(M_{E_i}, V_i)\in
G_{2e-f}^{e-f}(E_i)\}_{i=0}^{e-2a}$, which satisfies condition
(\ref{vanseq3}) with respect to $p_0$, by specifying the vanishing
sequences $a^{\mathfrak m_{E_i}}(p_i)$ for $0\leq i\leq e-2a$. For
$0\leq i\leq c-1$, the sequence $a^{\mathfrak m_{E_{i+1}}}(p_{i+1})$
is obtained from $a^{\mathfrak m_{E_{i}}}(p_{i})$ by raising all
entries by $1$, except for the term
$$a_{e-f-a}^{\mathfrak m_{E_{i+1}}}(p_{i+1})=a_{e-f-a}^{\mathfrak
m_{E_{i}}}(p_{i})=e-f-a+c.$$ After $c$ steps we arrive at the
following vanishing  sequence on $E_c$ with respect to $p_c$:
$$a^{\mathfrak m_{E_c}}(p_c)=(c, c+1, \ldots, e-f-a+c-1, e-f-a+c, r-a+2+c, r-a+3+c, \ldots, r+c+1).$$
For an index $c\leq i\leq e-2a-1$ which we write as $i=c+a\cdot
\beta +j$, with $0\leq j\leq a-1$ and $0\leq \beta\leq r-2-e+f$, we
choose $a^{\mathfrak m_{E_{i+1}}}(p_{i+1})$ to be obtained from
$a^{\mathfrak m_{E_{i}}}(p_{i})$ by raising all entries by $1$,
except for the term
$$a_{e-f-a+j+1}^{\mathfrak m_{E_{i+1}}}(p_{i+1})=a_{e-f-a+j+1}^{\mathfrak
m_{E_{i}}}(p_{i})=r-a+2+c+(a-1)\cdot \beta+2j.$$ In this way
$\mathfrak m\in \tilde{G}^{e-f}_{2e-f}(Y_0)$ becomes a (refined)
limit linear series which smooths to a complete linear series $M\in
G^{e-f}_{2e-f}(Y)$ on every smooth pointed curve $[Y, p]\in \cM_{e,
1}$ such that the ramification condition (\ref{vanseq3}) with
respect to $p$ is satisfied.

Next we construct a limit $\mathfrak g^1_{e-f+1}$ on $[Y_0, p_0]$,
say $\mathfrak b=\{(B_{E_i}, W_i)\in
G^1_{e-f+1}(E_i)\}_{i=0}^{e-2a}$ such that $a^{\mathfrak b}(p_0)=(0,
e-2f+1)$. For $0\leq i\leq e-2f$ we set $a^{\mathfrak
b_{E_i}}(p_i)=(i, e-2f+1)$. For an index of type $i=e-2f+2k-1$ where
$0\leq k\leq f-a$, we choose $a^{\mathfrak b_{E_i}}(p_i)=(e-2f+k-1,
e-2f+k+1)$. If $i=e-2f+2k$, we choose the sequence $a^{\mathfrak
b_{E_i}}(p_i)=(e-2f+k, e-2f+k+1)$. It is clear that each sequence
$a^{\mathfrak b_{E_i}}(p_i)$ is obtained from $a^{\mathfrak
b_{E_{i-1}}}(p_{i-1})$ by raising one entry by $1$ while keeping the
other fixed, hence $\mathfrak b$ is a limit $\mathfrak g^1_{e-f+1}$
which smooths to a pencil $B\in G^1_{e-f+1}(Y, p, (0, e-2f))$ on
every nearby smooth curve $[Y, p]$. For each $0\leq i\leq e-2a-1$,
there exists a section (unique up to scaling) $\sigma_i \in W_i$
such that
$\mbox{ord}_{p_i}(\sigma_i)+\mbox{ord}_{p_{i+1}}(\sigma_{i})=\mbox{deg}(B_{E_i})$.
We denote by $\sigma_i^c\in W_{i}$ a complementary section such that
$\{\mbox{ord}_{p_i}(\sigma_i),
\mbox{ord}_{p_i}(\sigma_i^c)\}=\{a_0^{\mathfrak b_{E_i}}(p_i),
a_1^{\mathfrak b_{E_i}}(p_i)\}$.

Using the set-up developed in \cite{EH3} and \cite{F2} for studying
degenerations of multiplication maps, we find that the assumption
that $\mu_{B, M}$ is not injective implies the existence elements
$0\neq \rho_i \in \mathrm{Ker}\{W_i\otimes V_i\rightarrow H^0(E_i,
B_{E_i}\otimes M_{E_i})\}$ for each $0\leq i\leq e-2a$, satisfying
the property that $\mbox{ord}_{p_{i+1}}(\rho_{i+1})\geq
\mbox{ord}_{p_{i}}(\rho_{i})+1$, for all $i$ (see e.g. \cite{F2}
Section 4, for an explanation of how to obtain the $\rho_i$'s).
Moreover, if
$\mbox{ord}_{p_{i+1}}(\rho_{i+1})=\mbox{ord}_{p_i}(\rho_i)+1$, then
if $\tau_i\in V_i$ is the section (unique up to scaling) such that
$\mbox{ord}_{p_i}(\tau_i)+\mbox{ord}_{p_{i+1}}(\tau_{i})=\mbox{deg}(M_{E_i})$,
then we must have that
$$\mbox{ord}_{p_i}(\rho_i)=\mbox{ord}_{p_i}(\tau_i)+\mbox{ord}_{p_i}(\sigma_i^c)=\mbox{ord}_{p_i}
(\sigma_i)+\mbox{ord}_{p_i}(\tau_i'),$$ where $\tau_i'\in V_i$ is
another section such that $\mbox{ord}_{p_i}(\tau_i')\neq
\mbox{ord}_{p_i}(\tau_i)$. In particular, since we have explicitly
described all the sequences $a^{\mathfrak b_{E_i}}(p_i)$ and
$a^{\mathfrak m_{E_i}}(p_i)$, the assumption that
$\mbox{ord}_{p_{i+1}}(\rho_{i+1})\leq \mbox{ord}_{p_i}(\rho_i)+1$
uniquely determines $\mbox{ord}_{p_i}(\rho_i)$.

Since $a^{\mathfrak b_{E_0}}(p_0)=(0, e-2f+1)$ and $\mu_{B_{E_0},
M_{E_0}}(\rho_0)=0$, the non-zero section $\rho_0$ must involve both
sections $\sigma_0$ and $\sigma_0^c$ and then clearly
$\mbox{ord}_{p_0}(\rho_0)\geq e-2f+1$. We prove inductively that for
all integers $0\leq i\leq e-2a$ we have the inequality
\begin{equation}\label{indu}
\mathrm{ord}_{p_i}(\rho_i)\geq e-2f+1+2i.
\end{equation}
Assuming (\ref{indu}) for $i\leq e-2a-1$, since
$\mbox{ord}_{p_{i+1}}(\rho_{i+1})\geq \mbox{ord}_{p_i}(\rho_i)+1$,
the only way (\ref{indu}) can fail for $i+1$ is when
$\mbox{ord}_{p_i}(\rho_i)=e-2f+2i+1$ and
$\mbox{ord}_{p_{i+1}}(\rho_{i+1})=\mbox{ord}_{p_i}(\rho_i)+1$. As
explained above, this implies that
$\mbox{ord}_{p_i}(\rho_i)=\mbox{ord}_{p_i}(\tau_i)+\mbox{ord}_{p_i}(\sigma_i^c)$.

Writing $i=c+a\cdot \beta +j$ as above, then
$\mbox{ord}_{p_i}(\tau_i)=r-a+2+c+(a-1)\cdot \beta +2j$\ if $i\geq
c$, while $\mbox{ord}_{p_i}(\tau_i)=e-f-a+c$, for $0\leq i\leq c-1$.
We deal only with the case $i\geq c$, the case $0\leq i\leq c-1$
being analogous. To determine $\mbox{ord}_{p_i}(\sigma_i^c)$ we must
distinguish between two cases: When $i=e-2f+2k-1$ with $k\geq 1$,
then $\mbox{ord}_{p_i}(\sigma_i^c)=e-2f+k-1$. Otherwise, we write
$i=e-2f+2k$ in which case $\mbox{ord}_{p_i}(\sigma_i^c)=e-2f+k+1$.
Suppose we are in the former case. Then we obtain the equality
$$e-2f+2i+1=\mbox{ord}_{p_i}(\rho_i)=\bigl(r-a+2+c+(a-1)\cdot
\beta+2j\bigr)+(e-2f+k-1),$$ which ultimately leads to the relation
$(a+2)(r-e+f-\beta)=a-j-1$. But $j\leq a-1$ and $\beta\leq r-e+f-1$,
hence we have reached a contradiction. The case when one can write
$i=e-2f+2k$ is dealt with similarly. All in all,  we may assume that
we have proved the inequality
$\mbox{ord}_{p_{e-2a}}(\rho_{e-2a})\geq e-2f+1+2(e-2a)$. We note
that on the curve $[T, q]=[E_{e-2a}, p_{e-2a}]$ we have that
$a^{\mathfrak b_{T}}(p_{e-2a})=(e-f-a, e-f-a+1)$, while
$$a^{\mathfrak m_{T}}(p_{e-2a})=(e-2a, e-2a+1, \ldots,
2e-f-3a, 2e-f-3a+3, \ldots, 2e-f-2a+2).$$  Equivalently $\mathfrak
b_T=|B|+(e-f-a)\cdot q$, where $B\in W_{a+1}^1(T)$, while $\mathfrak
m_{T}=(e-2a)\cdot q+|N|$, where $N\in \mbox{Pic}^{e-f+2a}(T)$ has
the property that $h^0\bigl(T, N(-(e-f-a+3)\cdot q)\bigr)\geq a$.
Remembering that $\mbox{ord}_{q}(\rho_{e-2a})\geq (e-2f+1)+2(e-2a)$,
after subtracting the base locus supported at $q$, we find an
element $$0\neq \rho_T\in \mbox{Ker}\{H^0(B)\otimes
H^0(N)\rightarrow H^0(B\otimes N)\}$$ such that
$\mbox{ord}_q(\rho_T)\geq e-f-a+1$. Equivalently, the multiplication
map $$\mu_{B, N}: H^0(B)\otimes H^0\bigl(N(-(e-f-a+3)\cdot
q)\bigr)\rightarrow H^0\bigl(B\otimes N(-(e-f-a+3)\cdot q)\bigr)$$
is not injective. By using Riemann-Roch we find that
$N(-(e-f-a+3)\cdot q)=K_T\otimes \tilde{B}^{\vee}$, where
$\tilde{B}\in W^1_{a+1}(T)$. Choosing $\tilde{B}=B\in W^1_{a+1}(T)$,
we notice that $\mu_{B, N}$ can be identified with the Petri map
$H^0(B)\otimes H^0(K_T\otimes B^{\vee})\rightarrow H^0(K_T)$ which
is injective because $[T]\in \cM_{2a}$ was chosen to be Petri
general. Thus we have reached a contradiction by reducing
(\ref{assumption2}) to the Gieseker-Petri theorem which completes
the proof in the case $(i)$.

Next we turn to case $(iii)$ when $[e/(r+1-e+f)]<2$. Since the
argument is similar to the one for $(i)$, we only outline the main
steps. If $e\leq r-e+f$, that is, when $\alpha_0=d-r-f-e$, we can
easily determine a linear series $(A, W_A)\in G^{r-e+f}_{d-e}(Y, p,
\overline{\alpha})$. Precisely, one can see that $A=K_Y\otimes
\OO_Y((d-3e+2)\cdot p)$ and $$|W_A|=(d-r-f-e)\cdot p+|K_Y\otimes
\OO_Y\bigl((r+f-2e+2)\cdot p\bigr)|.$$ In this case we have that
$|G_{d-e}^{r-e+f}(Y, p, \overline{\alpha})|=1$. Condition
(\ref{assumption2}) translates into saying that for a generic $(L,
W_L)\in G_{2e-f-1}^{e-f-1}(Y, p, \overline{\beta})$ we have the
vanishing statement
\begin{equation}\label{van3}
H^0\bigl(Y, L\otimes K_Y^{\vee}((e+f-2)\cdot
p)\bigr)=0\Leftrightarrow H^0\bigl(Y, K_Y^{\otimes 2}\otimes
L^{\vee}(-(e+f-2)\cdot p)\bigr)=0.
\end{equation}
One can prove (\ref{van3}) by degenerating $Y$ to a generic string
of elliptic curves and we skip the details. Finally, if
$[e/(r+1-e+f)]=1$, then $c=2e-r-f-1$ and condition
(\ref{assumption2}) boils down to showing that one can find a pencil
$B\in G^1_{e-c+1}\bigl(Y, p, (0, r-e+f-c+1)\bigr)$ and a linear
series $L\in G_{2e-f-1}^{e-f-1}(Y, p, \overline{\beta})$, such that
the multiplication map
$$H^0(B)\otimes H^0\bigr(K_Y^{\otimes 2}\otimes
L^{\vee}(-(2e-4-r)\cdot p)\bigl)\rightarrow H^0(K_Y^{\otimes
2}\otimes B\otimes L^{\vee}(-(2e-4-r)\cdot p)\bigr)$$ is injective.
This situation is handled along the lines of $(i)$ and we omit the
details. $\hfill$ $\Box$

Finally, we prove Theorem \ref{existence} assuming that condition (\ref{ass4})
is satisfied. This case is not covered by cases $(i)-(iii)$ above:

\begin{proposition}
Let $[Y, p]\in \cM_{e, 1}$ be a general pointed curve. Then there
exists a linear series $(A, W_A)\in G^{r-e+f}_{d-e}(Y, p,
\overline{\alpha})$ such that $H^0\bigl(Y, \OO_Y((d-1)\cdot p\otimes
A^{\vee})\bigr)=0$.
\end{proposition}

\begin{proof} By contradiction, we assume that $H^0(\OO_Y((d-1)\cdot p)\otimes A^{\vee})\neq 0$
for every $[Y, p]\in \cM_{e, 1}$ and for every linear series $(A,
V_A)\in G^{r-e+f}_{d-e}(Y, p, \overline{\alpha})$. We let $[Y, p]$
degenerate to the stable curve $[Y_0:=E_0\cup_{p_1}
E_1\cup_{p_2}\ldots \cup_{p_{e-3}}E_{e-3}\cup_{p_{e-2}}B, p_0]$,
where $E_0, \ldots, E_{e-3}$ are elliptic curves, the points $p_i,
p_{i+1}\in E_i$ are such that $p_i-p_{i+1}\in \mbox{Pic}^0(E_i)$ is
not a torsion class,
 and $[B, p_{e-2}]\in \cM_{2, 1}$ is such that $p_{e-2}\in B$ is
 not a Weierstrass point. For all integers
 $0\leq i\leq e-3$ we find that there exist sections
$$0\neq \tau_i\in H^0\bigl(\OO_{E_i}((d-1)\cdot p_i)\otimes A_{E_i}^{\vee}\bigr) \mbox{ and }0\neq
\tau_B=\tau_{e-2}\in H^0\bigl(\OO_{B}((d-1)\cdot p_{e-2})\otimes
A_B^{\vee}\bigr)$$ such that
$$0\leq \mbox{ord}_{p_0}(\tau_0)\leq
\mbox{ord}_{p_1}(\tau_1)\leq \ldots \leq
\mbox{ord}_{p_{e-3}}(\tau_{e-3})\leq
\mbox{ord}_{p_{e-2}}(\tau_{B}).$$ Moreover, we have that
$\mbox{ord}_{p_i}(\tau_i)\geq i$ for $0\leq i\leq e-2$. In
particular, $\mbox{ord}_{p_{e-2}}(\tau_{B})\geq e-2$. Since $\rho(e,
r-e+f, d-e, \overline{\alpha})=0$, limit $\mathfrak g^{r-e+f}_{d-e}$
on $E_0\cup \ldots \cup E_{e-3}\cup B$ are smoothable to every curve
of genus $g$. These finitely many limit $\mathfrak g^{r-e+f}_{d-e}$
are in bijective correspondence with possibilities of choosing the
vanishing sequences $\{a^{l_{E_i}}(p_i)\}_{0\leq i\leq e-3}$ and
$a^{l_B}(p_{e-2})$ in such a way that for all $0\leq i\leq e-3$, the
sequence $a^{l_{E_{i+1}}}(p_{i+1})$ is obtained from
$a^{l_{E_i}}(p_i)$ by raising all entries by $1$ except a single
entry which remains unchanged. To finish the proof it suffices to
exhibit a single limit $\mathfrak g_{d-e}^{r-e+f}$ on $E_0\cup
\ldots \cup E_{e-3}\cup B$ having the property that if $(A_B, V_B)$
denotes its $B$-aspect, then $H^0(\OO_B((d-e+1)\cdot p_{e-2})\otimes
A_B^{\vee})=0$.

We describe such a $\mathfrak g^{r-e+f}_{d-e}$ explicitly by
specifying the sequences $\{\alpha^{l_{E_i}}(p_i)\}_{0\leq i\leq
e-3}$ and $\alpha^{l_B}(p_{e-2})$. Clearly, $\alpha^{l_{E_0}}(p_0)$
equals $(\alpha_0, \ldots, \alpha_0,
\alpha^{l_{E_0}}_{r-e+f+1-c}(p_0)=\alpha_0+1, \ldots, \alpha_0+1)$.
For $1\leq i\leq c$, $\alpha^{l_{E_i}}(p_i)$ is obtained from
$\alpha^{l_{E_{i-1}}}(p_{i-1})$ by increasing all entries by $1$,
except for
$\alpha^{l_{E_i}}_{r-e+f+i-c}(p_i)=\alpha^{l_{E_{i-1}}}_{r-e+f+i-c}(p_{i-1})$.
Thus $\alpha^{l_{E_c}}(p_c)=(\alpha_0+c, \ldots, \alpha_0+c)$. Next,
for an index $i$ such that $c+\beta(r+1-e+f)< i\leq
c+(\beta+1)(r+1-e+f)$, where $0\leq \beta\leq [e/(r+1-e+f)]$, if we
write $i\equiv j+c \mbox{ mod } r+1-e+f$, with $1\leq j\leq r-e+f$,
the sequence $\alpha^{l_{E_i}}(p_i)$ is obtained from
$\alpha^{l_{E_{i-1}}}(p_{i-1})$ by raising all entries by $1$,
except for
$\alpha^{l_{E_i}}_{j-1}(p_i)=\alpha^{l_{E_{i-1}}}_{j-1}(p_{i-1})$.
Switching from ramification to vanishing sequences we obtain
$$a^{l_B}(p_{e-2})=(d-r-f-2, d-r-f-3, \ldots, d-e-5, d-e-4, d-e-2,
d-e-1),$$ that is, $A_B=\OO_B((d-e-2)\cdot p_{e-2})\otimes \mathfrak
g^1_2$, and then $$H^0(\OO_B((d-e+1)\cdot p_{e-2})\otimes
A_B^{\vee})=H^0(\OO_B(3\cdot p_{e-2})\otimes \mathfrak (\mathfrak
g^1_2)^{\vee})=0.$$ This contradicts the fact
$\mbox{ord}_{p_{e-2}}(\tau_B)\geq e-2$ which completes the proof.
\end{proof}

\section{Higher ramification points of a general line bundle}

In this section we prove Theorem \ref{genlb}. We fix an arbitrary
smooth curve $C$ of genus $g$ and for $n\geq 1$ we denote by
$[n]_C:\mbox{Pic}^d(C)\rightarrow \mbox{Pic}^{nd}(C)$ the
multiplication by $n$ map, $[n]_C(L):=L^{\otimes n}$. It is an
immediate consequence of Riemann-Roch that for a general $L\in
\mbox{Pic}^d(C)$, we have that $h^0(L^{\otimes
n})=\mbox{max}\{nd+1-g, 0\}$.

First we show that for a very general $L\in \mbox{Pic}^d(C)$ we have
that $w^{L^{\otimes n}}(p)\leq 1$ for all $p\in C$ and $n\geq 1$.
Indeed, let us assume that $w^{L^{\otimes n}}(p)\geq 2$, where $n$
is chosen such that $nd\geq g$, so that $h^0(C, L^{\otimes
n})=nd+1-g$. Then there are two possibilities:
 $$ \ (i)\   h^0\bigl(C, L^{\otimes
n}(-(nd+2-g)\cdot p)\bigr)\geq 1 \mbox{ or } \ (ii)\  \ h^0\bigl(C,
L^{\otimes n}(-(nd-g)\cdot p)\bigr)\geq 2.$$ In case $(i)$ we
consider the map $C\times C_{g-2}\rightarrow \mbox{Pic}^{nd}(C), \
(p, E)\mapsto \OO_C\bigl((nd+2-g)\cdot p+E\bigr)$ and we denote by
$\Sigma_n$ its image which is a divisor on $\mbox{Pic}^{nd}(C)$.
Then $(i)$ is equivalent to $L\in [n]_C^{*}(\Sigma_n)$ which is a
divisorial condition on $\mbox{Pic}^d(C)$ for each $n$.

In case $(ii)$ we look at the map $C \times C^1_g\rightarrow
\mbox{Pic}^{nd}(C), (p, E)\mapsto \OO_C\bigl((nd-g)\cdot p+E\bigr)$
and we denote by $V_n$ its image. Since $C^1_g$ is generically a
$\PP^1$-bundle over $C_{g-2}$, it follows that $V_n$ is a divisor on
$\mbox{Pic}^{nd}(C)$ and then possibility $(ii)$ is equivalent to
$L\in [n]_C^*(V_n)$. Thus we see that for $L\in
\mbox{Pic}^d(C)-\bigcup_{n\geq 1} [n]_C^*(\Sigma_n+V_n)$  all the
ramification points of all powers $L^{\otimes n}$ with $n\geq 1$,
are ordinary. This proves the first part of Theorem \ref{genlb}. To
prove the second part we start with the following:
\begin{proposition}
We fix a point $p\in C$ and integers $n$ and $d$ such that $nd\geq
g$. Then the locus
$$D_n:=\{L\in \mathrm{Pic}^d(C): h^0\bigl(C, L^{\otimes n}(-(nd+1-g)\cdot
p)\bigr)\geq 1\}$$ is an irreducible divisor on $\mathrm{Pic}^d(C)$
and $[D_n]=n^2 \theta$.
\end{proposition}
\begin{proof} We set $a:=\mbox{max}\{0, 2g-1-nd\}$ and define two
vector bundles $\E_n$ and $\F_n$ on $\mbox{Pic}^{d}(C)$ of the same
rank and having fibres $\E_n(L)=H^0(C, L^{\otimes n}\otimes
\OO_C(a\cdot p))$ and $\F_n(L)=H^0(C, L^{\otimes n}\otimes
\OO_{(a+nd+1-g)\cdot p}(a\cdot p))$ over each point $L\in
\mbox{Pic}^d(C)$. Then $D_n$ is the degeneracy locus of the morphism
$\E_n\rightarrow \F_n$ obtained by evaluation sections of
$L^{\otimes n}\otimes \OO_C(a\cdot p)$ along $(a+nd+1-g)\cdot p$.
The Picard bundle  $\E_n$ is negative (i.e. $\E_n^{\vee}$ is ample),
because $\E_n$ is the pull-back under the finite map $[n]_C$ of a
negative bundle on $\mbox{Pic}^d(C)$ (cf. \cite{ACGH}, pg. 310).
Moreover, $\F_n$ is algebraically equivalent to a trivial bundle,
hence $\E_n^{\vee}\otimes \F_n$ is ample too. Applying the
Fulton-Lazarsfeld connectedness theorem (see \cite{FL} or
\cite{ACGH} pg. 311), we conclude that $D_n$ is connected. Since
$D_n$ is also smooth in codimension $2$ we obtain that $D_n$ must be
irreducible. Finally, $[D_n]=c_1(\F_n-\E_n)=[n]_C^*(\theta)=n^2
\theta $.
\end{proof}

\noindent {\emph{End of the proof of Theorem \ref{genlb}}.} We fix
integers $1\leq a<b$ and consider the variety $\Sigma_{ab}:=\{(p,
L)\in C\times \mbox{Pic}^d(C): p\in R(L^{\otimes a})\cap
R(L^{\otimes b})\}$ and we denote by $\phi_1:\Sigma_{ab} \rightarrow
C$  and $\phi_2:\Sigma_{ab}\rightarrow \mbox{Pic}^d(C)$ the two
projections. For a fixed $p\in C$, the fibre $\phi_1^{-1}(p)$ is
identified with the intersection of the two irreducible divisors
$D_a$ and $D_b$. Since $[D_a]\neq [D_b]$ for $a\neq b$, it follows
that $D_a\cap D_b$ is of pure codimension $2$ inside
$\mbox{Pic}^d(C)$, therefore $\mbox{dim}(\Sigma_{ab})=g-1$. We
obtain that a line bundle $L\in \mbox{Pic}^d(C)-\bigcup_{a<b}
\phi_2(\Sigma_{ab})$ will enjoy the property that $R(L^{\otimes
a})\cap R(L^{\otimes b})=\emptyset$ for $a<b$. $\hfill$ $\Box$

\end{document}